\documentclass[12pt,twoside,leqno]{amsart}
\usepackage{amssymb,amsbsy,amsmath,amsfonts,amssymb,amscd,epsfig,
times,graphics,color,xy}

\definecolor{blue}{cmyk}{1.,1.,0.,0.67}
\definecolor{red}{cmyk}{0.,1.,1.,0.67}
\definecolor{green}{cmyk}{1.,0.,1.,0.67}
\usepackage[french]{babel}
\usepackage[T1]{fontenc} 
\newcommand{\C}{\mathbb{C}}

\newcommand{\N}{\mathbb{N}}
\newcommand{\R}{\mathbb{R}}

\newcommand{\blue}{\textcolor{blue}}
\newcommand{\green}{\textcolor{green}}
\newcommand{\red}{\textcolor{red}}

\makeatletter
\renewcommand{\@fnsymbol}[1]
{\ensuremath{\ifcase#1\or $*$\or $**$
\else\@ctrerr\fi}}
\makeatother

\begin{document}
\baselineskip=0.53cm

%\Large

$\:$\vspace{1cm}
\begin{center}
{\bf
INACTUALIT\'E ET INAD\'EQUATION

\medskip

DE LA PHILOSOPHIE DES MATH\'EMATIQUES 

\medskip

DE WITTGENSTEIN}

\bigskip

Jo\"el {\sc Merker}

\medskip

{\small\sf \'Ecole Normale Supérieure}

{\small\sf D\'epartement de Math\'ematiques et Applications}

{\footnotesize\bf www.dma.ens.fr/$\sim$merker/index.html}

\end{center}\medskip

\subsection*{ Prologue} 
Trop fréquemment, les mathématiques sont assimilées
lointainement et sans nuances\,\,---\,\,même par les philosophes les
plus ouvertement antiréalistes\,\,---\,\,à une entité immuable douée
d'une autonomie de principe par rapport au champ hétéronome de
l'expérience physique, biologique ou sociologique dans le monde. Mais
depuis plus d'une cinquantaine d'années, on ne sait plus comment
spéculer systématiquement sur le statut des théorèmes mathématiques en
tenant compte de manière globale et survolante de leur explosion et de
leur spécialisation, comme si ce monde dont la <<\,réalité\,>> encore
jugée problématique se <<\,réalisait\,>> sous nos yeux malgré
l'éternelle (et avantageuse) tentation de la dubitation philosophique,
qui s'est vue contrainte de se professionnaliser en se détachant des
sciences en action. Bien que le domaine philosophique contemporain
assigné comme tel semble en effet ne plus pouvoir suivre en pensée
cette complexification des contenus, et bien que la mathématique
engendre, de concert avec l'avancée du temps historique, ses propres
<<\,irréversibles-synthétiques\,>> dont le raffinement s'amplifie, il
est néanmoins du devoir de l'{\sl être-source du spéculatif
philosophique}\, de se confronter à l'{\sl être-ramifié du spéculatif
mathématique} pour y puiser des ressources dynamisantes et
structurantes. Si la mathématique pouvait enfin être reconnue comme
une philosophie dialectique réalisée et productrice à plein régime de
contenus argumentatifs authentiques, la philosophie générale devrait
lui emboîter le pas sans réticence et se défier méthodiquement des
cercles inactuels et reproductibles de raisonnements fermés.

Ce texte non théorisant n'a ici qu'une unique visée critique: limiter
la portée du général nominal dont on sait abuser lorsque les
raisonnements s'enlisent; centrifuger et éclater les cercles; ouvrir
et ramifier les questions en visant les niveaux contemporains; en un
mot: raffiner les analyses épistémologiques jusqu'au point où 
notre temps y éprouve ses résistances fertiles.

\subsection*{ Quel statut pour la proposition mathématique
non démontrée?} 
D'après la conception réaliste en mathématiques, la vérité des énoncés
mathématiques peut être réalisée sans que nous soyons en mesure de la
connaître, de la reconnaître ou de la démontrer pleinement: des
connexions rigoureuses existent toujours <<\,en réserve\,>>, dans un
<<\,matériau\,>> hypothético-déductif <<\,potentiellement indéfini\,>>
que l'on <<\,explore\,>> ou que l'on <<\,découvre\,>>. Tout \`a
l'opposé des visions <<\,réalistes\,>> et pour se démarquer des
<<\,na\"{\i}vetés\,>> qui les accompagnent, Wittgenstein affirme à
divers endroits de ses remarques philosophiques qu'il existe une
différence de nature profonde entre les propositions
<<\,pressenties\,>>, non démontrées, et celles qui sont déjà insérées
dans une grammaire formelle et autonome, conçue comme système
prédéfini de règles du jeu tel que l'arithmétique, ou la théorie des
ensembles ou la topologie générale. Ce fossé conceptuel majeur est
incontournable, inexpugnable et ne peut en aucune façon être comblé.

\smallskip
\hfill
\begin{minipage}[t]{11cm}
\baselineskip=0.37cm
{\small{\sf\blue{
La démonstration se distingue radicalement de la vérification d'une
proposition ordinaire déjà comprise, en ce sens que <<\,la
démonstration fait partie de la grammaire de la proposition\,>> (PG,
p.~370); de sorte que <<\,la proposition avec sa démonstration
appartient à une tout autre catégorie que la proposition sans la
démonstration\,>> (PG, p.~371). La proposition non démontrée ne
représente pas un fait mathématique dont nous ne savons pas encore et
dont nous cherchons à savoir s'il est réalisé ou non. Son sens est
uniquement celui d'une incitation à la recherche mathématique et d'une
directive ou d'une suggestion pour la recherche. Mais, alors qu'une
hypothèse empirique (par exemple, une hypothèse médicale) conserve le
même sens, lorsqu'elle est vérifiée, la démonstration mathématique
modifie la position de la proposition dans le langage lui-même, et
donc son sens.
}}
\hfill
\green{\cite{ bo1988},~p.~60.}
}
\end{minipage}\medskip

Différence de catégorie parce que seul l'effectivement démontré est
inséré dans l'architecture syntaxique des énoncés mathématiques. Mais
rien de tel pour le non-démontré: il est suspendu dans l'hypothétique,
indécis et non sanctionné; il est <<\,autre\,>>, parce qu'il
appartient à un autre univers de langage et de pensée.
Voilà une démarcation bien claire et bien nette.

\subsection*{ L'irréversible-synthétique en mathématiques} 
Pour l'intuition de compréhension, aucune difficulté à résumer cette
thèse: les oppositions de nature sur lesquelles on insiste en
philosophie appartiennent en effet aux informations les plus
immédiatement saisissables par la pensée <<\,archaïque et
reptilienne\,>> qui contrôle toutes les mobilisation neuronales du
cerveau, donc de la pensée humaine effective. Répéter un {\it
distinguo}\, est non seulement autorisé, mais cela est aussi
nécessaire; c'est une manière de chercher à circonscrire un fait
essentiel, en le soumettant, sans l'épuiser, à une méditation variée
qui éclaire les premières parois d'un fossé conceptuel.

\smallskip

Réexprimons donc l'idée en d'autres termes. Si l'on admet sans
discuter, comme
semble le faire Wittgenstein, que l'univers du non syntaxiquement
sanctionné ne doit pas faire l'objet d'une recherche visant à lui
conférer structures internes et logiques autonomes qui entretiennent
des liens complexes et délicats avec l'ordre du discours formalisé,
une chose est pour l'instant certaine: le champ mathématique est
universellement traversé par ce que nous appellerons dorénavant
l'\green{\sl irréversible-synthétique}. Est {\sl synthétique}\, tout
raisonnement qui compose avec des objets de pensée et qui rassemble
des éléments de connaissance en un tout cohérent, en travaillant de
manière locale ou (partiellement) globale. Est {\sl irréversible}\,
tout phénomène physique qui ne fonctionne que dans un seul sens, sans
pouvoir être renversé spontanément, comme par exemple la formation
d'un précipité chimique ou l'oxydation du fer. Mais dans le domaine
abstrait, l'irréversible ne peut pas être caractérisé en termes
organiques, ou être quantifié en termes d'entropie. Parler
d'<<\,irréversibilité mathématique\,>> ne constituerait certainement
pas une expression adéquate, parce qu'il n'y a pas, dans le domaine de
la pensée pure, de démonstrations en marche par elles-mêmes qu'il
suffirait de déclencher en confrontant les définitions aux questions,
dans un creuset magique et hypothétique que personne n'a encore
découvert. 

\medskip

\hspace{-2.75cm}
\fbox{\red{\small\sf En mathématiques, nul automatisme empirique, et il
n'y a pas d'essence motrice séparée.}}

\medskip

Aussi l'<<\,irréversible\,>>
doit-il se rapporter dans sa notion propre\footnote{\, Le biologique
du rationnel doit notamment y consacrer son potentiel
d'irréversibilité.} à ce qui fait que l'essence des mathématiques est
de démontrer synthétiquement, chaque démonstration {\it synthétique}\,
faisant {\it basculer les contenus de manière irréversible dans le
champ expansif des résultats rigoureusement établis}.

\smallskip

En résumé, la thèse forte autour de laquelle se concentre la pensée de
Wittgenstein dit simplement que l'irréversible-synthétique domine le
statut de la proposition mathématique: c'est dans l'{\it a
posteriori}\, d'une synthèse, et seulement dans cet {\it a
posteriori}, que s'affirme la signification mathématique d'un énoncé.

\subsection*{ Insuffisances spéculatives}
Toutefois, un réel danger de circularité menace toute position qui
affirme unilatéralement une thèse d'opposition, quelle qu'en soit la
portée. Parce qu'elle s'inscrit dans la temporalité propre du monde,
l'opposition fondamentale entre l'<<\,avant\,>> et l'<<\,après\,>>
appartient en effet aux dialectiques les plus évidentes et les plus
omniprésentes pour ce qui est de la vie continue de l'esprit. D'un
point vue spéculatif, on ne peut pas se cantonner à répéter cette
constatation chaque fois qu'on la voit se manifester dans la vie
propre des étants abstraits et concrets que l'on fréquente, {\it parce
que}\, de très nombreuses questions théoriques invitent à explorer les
failles imprécises de ce <<\,fossé conceptuel\,>> qui sépare
démonstrations achevées et supputations provisoires.

\medskip 

$\bullet$
Quand peut-on parler de validation définitive d'un résultat
mathématique? Quel critère choisir? Quelle ligne de démarcation
proposer?

\medskip 

$\bullet$
La pensée du conjectural doit-elle être considérée par principe comme
définitivement éliminée à l'instant même où la proposition
mathématique formalisée confirme l'attente de vérité? Comment alors
s'effectue une telle <<\,cristallisation-élimination\,>>?

\medskip 

$\bullet$
\`A quel moment peut-on être certain que la proposition s'incrit
véritablement dans le système grammatical autorisé? Doit-on établir
des nuances en fonction de la structure, de la longueur et de la
complexité des preuves? Si l'on décompose un théorème donné en
fractions partiellement ou totalement vérifiées, doit-on être conduit
à parler d'hétéronomie du champ démonstratif?

\medskip 

$\bullet$
Lorsqu'il est soumis à révision (correction), comment un théorème donné
modifie-t-il son inscription dans la grammaire générale des énoncés
mathématiques?

\medskip 

$\bullet$
Quel statut donner aux preuves formelles qui ont été publiées dans des
revues de mathématiques internationales, mais qui se sont en vérité
avérées incorrectes après examen ultérieur, et souvent imprévisible,
par d'autres mathématiciens? Le philosophe du <<\,fossé conceptuel\,>>
a-t-il été victime d'une illusion, d'un mirage? Quand et comment
peut-il être certain qu'il s'en rend compte\footnote{\, Wittgenstein
dit bien: <<\,\blue{\sf La proposition avec sa démonstration
appartient à une tout autre catégorie que la proposition sans la
démonstration}\,>>. Onze années séparent les deux <<\,preuves\,>> du
théorème des quatre couleurs publiées par Kempe en 1879 et Tait en
1880 de la découverte par Petersen en 1891 de la présence d'un
<<\,trou\,>> tellement important qu'il fallut encore attendre 1976
(Appel et Haken, après des idées décisives de Birkhoff, Heesch et
d'autres) pour que l'on domine la <<\,zoologie\,>> des milliers de
noyaux <<\,inévitables\,>> qui remplaçaient le <<\,centre\,>> de la
carte dont était parti Kempe. }?

\medskip 

$\bullet$
\`A quel moment <<\,{\it cela}\,>> bascule-t-il? et à quel
moment <<\,{\it cela}\,>> re-pivote-t-il en arrière en cas d'erreur?
Où et quand mémoriser l'erreur? Quel statut lui réserver?

\subsection*{ Dérobade philosophique}
Or face à de telles questions préliminaires, Wittgenstein semble
choisir de se soustraire sciemment, intentionnellement au devoir
d'analyser et de penser la complexité des liens qui unissent la pensée
intuitive, prospective et informelle au régime d'appropriation réglée
qu'offrent axiomatisation et formalisation.

\renewcommand{\thefootnote}{\fnsymbol{footnote}}

\smallskip
\hfill
\begin{minipage}[t]{11cm}
\baselineskip=0.37cm
{\small{\sf\blue{
La proposition mathématique non démontrée ne contient pas une
anticipation d'un fait qui a pu être suggéré par des expériences et
dont la démonstration se chargera d'établir l'existence. Wittgenstein
dit d'elle qu'elle est <<\,un poteau indicateur pour la recherche
mathématique, une incitation à des constructions mathématiques\,>>
(PG, p.~371). Ce qui lui donne pour l'instant un sens mathématique est
essentiellement le complexe de résonances, d'associations,
d'analogies, {\it etc.} qu'elle suscite dans le système des
mathématiques et qui fournit à la fois un stimulant et une direction à
la recherche.
}}
\hfill
\green{\cite{ bo1988},~p.~195.}
}
\end{minipage}\medskip

\renewcommand{\thefootnote}{\arabic{footnote}}

Poteau indicateur, poteau télégraphique, poteau-frontière, résonances
imprécises, analogies\,\,---\,\,voilà ce à quoi semble être réduite la
pensée en acte dans la recherche effective. Boîte noire, dirons-nous
tout simplement: opaque au philosophe, elle fonctionne à une distance
éloignée de lui; et c'est bien de la tête <<\,noire\,>> du
mathématicien qu'il est question ici; pourtant, la tâche que s'assigne
le philosophe wittgensteinien met {\it a priori}\, entre parenthèses
le devoir de comprendre et de fréquenter les réseaux de raisonnements
qui {\it produisent}\, les constructions mathématiques
formalisées. C'est cela que pourrait être tenté de lui reprocher tout
mathématicien professionnel habitué à {\it jongler entre les deux
niveaux en diluant des frontières imprécises}, habitué à métamorphoser
des briques de rigueur formelle en {\it forces argumentatives douées
de mobilité questionnante}, habitué à {\it vivre pendant des années en
compagnie de problèmes ouverts partiellement explorés}, et dont il est
incapable, tout autant que ses collègues et concurrents directs,
d'évaluer l'horizon de difficulté rémanente, c'est bien cela en effet
qu'on serait tenté d'opposer à Wittgenstein, si son discours portait
véritablement sur les mathématiques tout entières, comme semblerait le
prétendre son vocabulaire qu'il ne veut jamais spécifique ou
spécialisé; et même sans chercher la polémique\,\,---\,\,je dirais
même plus, sans chercher l'affrontement avec des adversaires peut-être
inconsistants à qui les vraies difficultés, tant qu'ils s'y dérobent,
restent invisibles\,\,---\,\,on pourrait de surcroît se demander
véritablement comment il a pu être possible, dans l'histoire des
idées, que la parole très assertorique de certains qui n'ont jamais
<<\,créé\,>> de mathématiques soit parvenue à énoncer et à faire
circuler un discours qui puisse faire autorité, dans certains milieux
philosophiques, quant à la manière dont les mathématiques <<\,se\,>>
créent, l'emploi du pronominal réfléchi <<\,se\,>> montrant ici {\it
cum grano salis}\, à quel point il s'agit d'une <<\,boîte\,>>
absolument <<\,noire\,>> pour ceux qui tentent d'échafauder un
discours universel à ce sujet.

\subsection*{ Liberté mathématique}
Aucun discours universel sur les actes de
<<\,basculement\,>> vers des résultats vrais, nouveaux et sanctionnés
n'est parvenu à s'ériger, à aucune période de l'histoire des
mathématiques, en tant que champ de principes stables, série de
méthodes directives, ou ordre réglé de découvertes potentielles. 

\smallskip
\hfill
\begin{minipage}[t]{11cm}
\baselineskip=0.37cm
{\small{\sf\blue{
Les mathématiques sont un outil de liberté.
}}
\hfill
\green{Adrien~{\sc Douady}.}
}
\end{minipage}\medskip

\noindent
Tous les mathématiciens professionnels savent la mathématique trop
libre de par son caractère imprévisible, notamment parce qu'ils
éprouvent journellement l'<<\,exoticité\,>> et
l'<<\,incompréhensibilité\,>> de tous les résultats mathématiques qui
sont éloignés de ce qu'ils connaissent de très près. Naïvetés, donc,
que les séduisantes formules
wittgensteiniennes!\,\,---\,\,lorsqu'elles extrapolent dogmatiquement
leur portée en abusant de généralité terminologique.

Au contact des mathématiques contemporaines et s'il se décidait à les
fréquenter véritablement et à les analyser dans leur complexité
actuelle, Wittgenstein ressuscité démultiplierait peut-être sa pensée,
mais on ne pourrait pas alors tout à fait exclure que ses lecteurs
épigones ne puissent plus être à même d'étudier ses travaux pour 
asseoir une autorité philosophique. Imaginons en effet ses
paragraphes compacts de deux à vingt ou trente lignes\,\,---\,\,si
faciles à lire pour le commun des philosophes\,\,---\,\,se
métamorphoser en centaines de pages ciselées qui exigent des années de
formation mathématique préalable?

\medskip

\hspace{-2.05cm}
\fbox{\fbox{\small\sf\red{ L'irréversible mathématique
doit forcer à complexifier les règles du jeu exégétique.}}}

\subsection*{ Conjectures expérimentales étrangères aux démonstrations
rigoureuses} Après cette brève contre-argumentation,
reprenons l'examen des thèses wittgensteiniennes au sujet de
l'induction en mathématiques.

\smallskip
\hfill
\begin{minipage}[t]{11cm}
\baselineskip=0.37cm
{\small{\sf\blue{\label{fosse}
Wittgenstein soutient qu'il existe un gouffre conceptuel
infranchissable entre la conjecture, qui anticipe les résultats d'une
série d'expériences de calcul hypothétiques, et la démonstration, qui
prescrit, de façon complètement impersonnelle et intemporelle, quelque
chose à propos des résultats en question. La première, pour autant
qu'elle ressemble à ce qu'on appelle ordinairement une conjecture, dit
simplement qu'aucun contre-exemple ne se présentera, la seconde exclut
que quelque chose puisse être appelé un contre-exemple.
}}
\hfill
\green{\cite{ bo1988},~p.~194.}
}
\end{minipage}\medskip

Effectivement, la différence est radicale: rappelons par exemple le
destin <<\,attentiste sur plus d'un siècle\,>> de la loi {\it
quantitative}\, de répartition des nombres premiers\footnote{\, {\it
cf.} {\it e.g.} J.-P. {\sc Delahaye}, {\em Merveilleux nombres
premiers. Voyage au c{\oe}ur de l'arithmétique}. Belin, Paris,
2000. }: par des arguments élémentaires, Legendre a montré en 1808 que
l'ensemble des nombres premiers admet une densité nulle sur $\N = \{
1, 2, 3, \dots \}$; mais comme Euler avait déjà établi auparavant que
la somme des inverses des nombres premiers diverge: $\sum_{ p \,
\text{\rm premier}}\, \, \frac{ 1}{ p} = + \infty$, cette densité
nulle ne pouvait pas signifier une très forte
raréfaction. Existe-t-il alors une loi
mathématique qui décrit cette
raréfaction de manière quantitative?

Ce fut semble-t-il dès 1792 qu'à l'âge de 15 ans, Gauss émit la toute
première {\it hypothèse quantitative précise}\, de
raréfaction\footnote{\, Ce fait est attesté en 1848 dans une 
réponse de Gauss à l'astronome allemand Johan Encke 
qui aurait découvert une
loi similaire; les mentions éparses que Gauss
formulaient dans sa maturité quant à ses découvertes de jeunesse sont
à prendre très au sérieux, étant donné qu'il se refusait à publier la
plupart de ses résultats partiels, et {\it a fortiori}\, les
conjectures qu'il n'était pas parvenu à démontrer.}: en examinant les
tranches de 1\,000 entiers dans les tables de nombres premiers (qu'il
corrigeait au passage jusqu'à des entiers dépassant plusieurs
millions), Gauss observa qu'au voisinage d'un entier $n$ quelconque,
la densité des nombres premiers est de l'ordre de $\frac{ 1}{ {\rm
log}\, n}$. Alors il émit l'hypothèse que le nombre $\pi ( n)$ de
nombres premiers inférieurs ou égaux à $n$ devrait être
asymptotiquement égal au logarithme intégral\footnote{\, Extraits de
la lettre de Gauss à Encke, 24 décembre 1849, traduite en anglais
dans: L.J.~{\sc Goldstein}, {\em A history of the prime number
theorem}, Amer. Math. Monthly {\bf 80} (1973), 599--615: <<\,Your
remarks concerning the frequency of primes were of interest to me in
more ways than one. You have reminded me of my own endeavors in this
field which began in the very distant past, in 1792 or 1793 after I
have acquired the Lambert supplements to the logarithmic tables.
[\dots] I counted the primes in several chiliads [\dots]. I soon
recognized that behind all of its fluctuations, this frequency is on
the average inversely proportional to the logarithm, so that the
number of primes below a given bound $n$ is approximately equal to
$\int \, \frac{ dn}{ {\rm log}\, n}$, where the logarithm is
understood to be hyperbolic. Later on, when I became acquainted with
the list in Vega's tables (1796) going up to 400\,031, I extended my
computations further, confirming that estimate. In 1811, the
appearance of Chernau's cribrum gave me much pleasure and I have
frequently (since I lack the patience for a continuous count) spent an
idle quarter of an hour to count another chiliad here and there;
although I eventually gave it up without quite getting through a
million. Only some time later did I make use of the diligence of
Goldschmidt to fill some of the remaining gaps in the first million
and to continue the computation according to Burkhardt's tables. Thus
(for many years now) the first three millions have been counted and
checked against the integral.

\medskip

\begin{center}
\begin{tabular}[t]{cccccc}
$n$ & $\pi(n)$ & $\int\,\frac{dn}{{\rm log}\,n}$ & {\sc Error}
& {\sc Your formula} & {\sc Error}
\\
\hline
\\
\vspace{-7mm}
\\
500\,000 & 41\,556 & 41\,606,4 & +50,4 & 41\,596,9 & +40,9
\\
1\,000\,000 & 78\,501 & 79\,627,5 & +126,5 & 78\,672,7 & +171,7
\\
1\,500\,000 & 114\,112 & 114\,263,1 & +151,1 & 114\,374,0 & +264,0
\\
2\,000\,000 & 148\,883 & 149\,054,8 & + 171,8 & 149\,233,0 & +350,0
\\
2\,500\,000 & 183\,016 & 183\,245,0 & +229,0 & 183\,495,1 & 479,1
\\
3\,000\,000 & 216\,745 & 216\,970,6 & +225,6 & 217\,308,5 & +563,5
\end{tabular}
\end{center}

\medskip\noindent
[\dots] The chiliad from 101\,000\,\,---\,\,102\,000 in Lambert's
Supplement is virtually crawling with errors; in my copy, I have
indicated seven numbers which are not primes at all, and supplied two
missing ones. [\dots].\,>> } $\int_2^n \, \frac{ dt}{ {\rm log}\, t}$.
L'approximation équivalente un peu moins précise $\pi ( n) \sim \frac{
n}{ {\rm log} \, n}$ a été conjecturée\footnote{\, En fait, dès 1798,
Legendre affirmait que l'on a exactement $ \pi(n) = \frac{n}{{\rm
log}\,n+A(n)}$, <<\,où $A ( n)$ est approximativement égal à $1,
08366\cdots$\,>>. Mais cet énoncé incorrect devait être mis en défaut
assez rapidement. }  par Legendre en 1808.

Premier moment expérimental, donc, purement observationnel et
simplement cantonné à un suivi comptable; patience obstinée de
calculateur prodige et ingénu était ici requise\footnote{\, Gauss a
donc poursuivi ces recherches bien des années après avoir publié ses
{\it Disquisitiones Arithmeticae}. Mille pages environ sont
nécessaires pour écrire ces 216\,745 nombres premiers \`a raison de
deux cent dix nombres premiers par page sur trois colonnes. Compter $\pi (
n + 1000) - \pi ( n)$ est immédiat. Calculer une valeur numérique
précise du logarithme intégral prend quelque temps. S'assurer que les
tables ne comportent pas d'erreur est beaucoup plus délicat. }. Pour
tester ou deviner des lois plausibles, il est {\it nécessaire}, sinon
incontournable, d'{\it \'eriger au préalable d'arides pyramides
numériques}\, pour en scruter les structures translucides noyées dans
une opacité primordiale\footnote{\, Voici un autre exemple célèbre que
le labyrinthe de l'induction nous a transmis dans l'histoire des
mathématiques.  Soit $n$ un entier naturel $\geqslant 1$. Question: de
combien de manières distinctes peut-on casser $n$ en morceaux
discrets, {\it i.e.} écrire $n = a + b + c + \cdots$, où $a$, $b$,
$c$, $\dots$ sont des entiers $\geqslant 1$? En fait, il y a {\it
deux}\, questions, suivant que l'on décide (ou non) de tenir compte de
l'ordre dans lequel sont écrits les constituants de $n$. Avec
distinction de l'ordre de sommation, la réponse est essentiellement
{\it trop}\, simple: une démonstration par récurrence montre en effet
qu'il y a juste $2^{ n-1}$ possibilités, par exemple: $3 = 2 + 1 = 1 +
2 = 1 + 1 + 1$. Mais les choses sont incroyablement plus compliquées
lorsqu'on néglige l'ordre; notons donc $p ( n)$ le nombre de
\green{\sl partitions}\, de $n$; par exemple pour $n = 5$, on a $p (
5) = 7$, car $5 = 4 + 1 = 3 + 2 = 3 + 1 + 1 = 2 + 2 + 1 = 2 + 1 + 1 +
1 = 1 + 1 + 1 + 1 + 1$.

\smallskip
\hfill
\begin{minipage}[t]{11cm}
\baselineskip=0.37cm
\blue{\scriptsize\sf
There is a famous story concerning the search for some kind of pattern
in the table of the $p ( n)$'s. This is told of Major Mac Mahon who
kept a list of these partition numbers arranged one under another up
into the hundreds. It suddenly occured to him that, viewed from a
distance, the outline of the digits seemed to form a parabola! Thus
the number of digits in $p ( n)$, the number of partitions of $n$, is
around $C \sqrt{ n}$, or $p ( n)$ itself is very roughly $e^{ \alpha
\sqrt{ n}}$. The first crude assessment of $p ( n)$!}

\blue{\scriptsize\sf 
Among other things, however, this does not tell us not
to expect any simple answers. Indeed later research showed that the
true asymptotic formula for $p ( n)$ is $\frac{ e^{ \pi \sqrt{ 2n /
3}}}{ 4 \sqrt{ 3}\, n}$!
}
\hfill
{\small\green{D.J.~{\sc Newman}.}}
\end{minipage}}\,\,---\,\,sinon, quelle vision transcendante
viendrait secourir l'intuition prospective? Et actuellement, la
théorie dite <<\,computationnelle\,>> des nombres regorge de
conjectures observationnelles quantitatives parfaitement certaines,
sans qu'aucune des expériences numériques automatisées lancées sur des
ordinateurs super-performants ne puisse offrir d'indication quant à un
hypothétique champ démonstratif afférent: raison est donc donnée à
Wittgenstein sur ce point, si l'on s'en tient aux exemples pour
lesquels l'inconnu déductif reste sensiblement à l'écart du scruté
expérimental.

\subsection*{ Inexactitudes et expressions inappropriées}
Toutefois, dans le court extrait de \green{ \cite{ bo1988}} reproduit
ci-dessus p.~\pageref{fosse}, la manière même de s'exprimer est
imparfaite et inadéquate, voire tout simplement <<\,
\underline{\red{\bf fausse}}\,>>, si l'on doit s'autoriser à employer,
au sein d'un débat de philosophie des mathématiques, une terminologie
typique de la pratique des mathématiciens.

Tout d'abord, l'adjectif <<\,\blue{\small\sf infranchissable}\,>> dans
l'expression <<\,\blue{\small\sf gouffre conceptuel
infranchissable}\,>> est absurde: au contraire, certaines conjectures
ont été, sont et seront démontrées. Justement les mathématiciens
inventent des concepts dont ils <<\,remplissent\,>> ces <<\,fossés
conjecturaux\,>> jusqu'à pouvoir les {\it franchir}. Toute la
difficulté est de pouvoir penser ce mouvement complexe et mystérieux. Place
aux paradoxes, aux questions et à la philosophie!

\smallskip
\hfill
\begin{minipage}[t]{11cm}
\baselineskip=0.37cm
{\small{\sf\blue{
The brain of every mathematician carries a fragment of our ``cloud in
the tree'', a little personal cloud where our synapses touch Hilbert's
tree. These little clouds may have fractal geometry and thus
relatively large boundaries. Hard to tell at this stage but one may
take analogy from the study of the human movements where our neuron's
system {\it avoids}\, most paths through many degrees of freedom as
experiments show. This may be also the mathematical strategy of our
brain, responsible for instance, for the equality $P = NP$ of everyday
mathematics. We solve our problems essentially as fast as we state
them. It took, probably, a couple of thousand brain-hours to state
the Fermat theorem and mere instance (compared to $\exp 2000$) to
solve it, no more than $10^5$ brain-hours. (Actually, one has to
compare the length of the proof to the time needed to find it. Maybe,
the {\it shortest}\, proof of Fermat in reasonable units is of the
order ${\rm log}$(time spent on the search of the proof).) This
``practical equality'' $P = NP$ is in flagrant contradiction with our
mathematical intuition as we expect $NP$ to be far away from $P$.
Here is a fundamental gap in our understanding (if there any) of how
mathematics works. We need, besides pure thought, biological,
psychological study and/or computer experimentation. But as a
community we shy away from such problems, scared of contamination by
philosophy.
}}
\hfill
\green{Mikhail {\sc Gromov}.}
}
\end{minipage}\medskip

Ensuite, la démonstration mathématique ne <<\,\blue{\small\sf
prescrit}\,>> pas\footnote{\, Nous poursuivons l'analyse critique de
l'extrait en question, p.~\pageref{fosse}. }, elle {\it établit} (des
propositions, des résultats, des théorèmes). On pourra certes admettre
qu'elle <<\,empile\,>> des arguments, qu'elle <<\,aligne\,>> des
raisonnements, qu'elle <<\,combine\,>> des techniques diverses. Mais
le terme inapproprié <<\,prescrire\,>> est vraiment à proscrire, ne
serait-ce que parce qu'il suggère quelque chose comme une décision de
législateur ou un acte médical, bref une espèce de recommandation
expresse, d'exigence, d'obligation ou d'ordre qu'il serait
déraisonnable et fou de proférer face à certains problèmes très
ouverts qui se posent avec tout un réservoir de potentialités
imprévisibles.

Même s'il doit s'agir de règles, d'axiomes, ou de déclarations
syntaxiques de concepts, <<\,prescrire\,>> ne peut en aucun cas
absorber la portion principale de l'énergie de recherche en
mathématiques. En effet, le champ mathématique n'est pas simplement
<<\,prescrit\,>> ou <<\,déclaré\,>> par des démonstrations ou par des
règles de langage, fussent-elle dûment établies avec toute la rigueur
formelle, parce que la <<\,prescription\,>> ou plutôt 
la <<\,déclaration\,>> et la <<\,mise en place\,>>
des <<\,règles\,>> n'est qu'un prologue au déploiement du champ de
l'irréversible-synthétique, qui s'éclaire ensuite grâce à des quanta
argumentatifs articulés et mobilisés dynamiquement {\it dans}\, les
démonstrations.

De plus, affirmer que la <<\,\blue{\small\sf démonstration
mathématique prescrit [\dots] quelque chose à propos des résultats en
question}\,>> constitue une périphrase raccourcie, maladroite et trop
rapide pour nommer le lien complexe qui unit les énoncés aux
arguments, comme si ce que la démonstration dévoile du résultat
qu'elle démontre devait subir un dédoublement et se constituer en même
temps comme une réalisation exemplaire du point de vue grammatical
prescriptif; comme si les chaînes formalisées d'arguments lançaient,
dans le champ de l'indéfini axiomatisé, un éclair qui se pétrifierait
du même coup pour confirmer l'immanence fixée de l'univers des règles; bref,
comme si toute démonstration mathématique devait nécessairement être
entraînée dans une métaphysique wittgensteinienne.

Ici encore, la généralité du vocabulaire invite à négliger la
complexité des situations: le caractère toujours partiel des saisies
axiomatico-formalistes dans la pratique mathématique, et la permanence
des horizons imprécis de questionnement font que la démarcation même
entre l'argumentatif et le déductif se fragmente à la fois dans
l'histoire d'une spécialité et dans l'appréhension mentale des
théorèmes. En tout cas, que l'on n'objecte pas que la rapidité
d'exécution des phrases examinées expose inévitablement à certaines
imperfections, car il s'agit bien ici d'un des problèmes les plus
difficiles de la philosophie des mathématiques: penser la réalisation
de l'irréversible-synthétique. Finesse de la spéculation et précision
dans la terminologie doivent être d'emblée exigées.

Continuons: que fait la conjecture? Non, elle n'<<\,\blue{\small\sf
anticipe}\,>> pas <<\,\blue{\small\sf les résultats d'une série
d'expérience de calculs hypothétiques}\,>>! Même en se restreignant
aux aspects purement expérimentaux de la théorie des nombres (le
conjectural s'exerce en fait dans toutes les spécialités
mathématiques), il serait fort réducteur de n'y voir qu'une prévision
tout expérimentale d'expériences numériques futures. Bien que la
phrase citée soit contrainte ici de continuer à maintenir une nette
démarcation afin de garantir la cohérence locale de sa thèse, il nous
faut rappeler que la conjecture énonce des règles, prétend des
régularités, soupçonne des théorèmes, devine des lois, et s'exprime la
plupart du temps dans le même langage formalisable que toutes les
propositions qui sont dûment établies dans le sanctuaire
hypothético-déductif. Une conjecture ordinaire, c'est un énoncé sans
démonstration, l'énoncé vraiment possible d'un théorème vraisemblable
que l'on pose dans un moment de suspens face à de l'inconnu qui
résiste. Par définition, la conjecture est un énoncé potentiel fort
d'une pensée structurée, bien que transversale au régime réglé des
grammaires formelles, c'est un énoncé qui appelle une démonstration,
ou qui subira une réfutation.

Dans la communauté internationale des mathématiciens, rares sont les
conjectures qui s'affirment comme citadelles de pensée résistant à de
multiples assauts intellectuels: en un mot, rares sont les conjectures
dignes de ce nom, parce que la conjecture requiert d'embrasser des
abysses synthétiques spécifiques qui focalisent un fort enjeu
mathématique et face auxquelles on doit se sentir écartelé et trop
faible pour s'autoriser à en dire quelque chose.

\smallskip
\hfill
\begin{minipage}[t]{11cm}
\baselineskip=0.37cm
\blue{\small\sf
As you said, Don [Zagier], the conjecture is the most responsible
thing one can do and sometimes people make conjectures when they
absolutely have no right to make conjectures. A conjecture really
comes hard. I agree with you, Don, that one could make a serious
conjecture once or twice in one's life, after deep thinking. You come
to a deep understanding, and you cannot finish it, and you make a
conjecture. You just cannot turn any question into a conjecture.
}
\hfill
{\small\green{Mikhail {\sc Gromov}}}
\end{minipage}\medskip

Il est par ailleurs surprenant de lire dans le
même extrait (p.~\pageref{fosse} {\it supra}) que la
conjecture <<\,\blue{\small\sf dit simplement qu'aucun contre-exemple
ne se présentera}\,>>, car dans la forme même de son énonciation, la
conjecture ne s'attarde en général pas à se contraposer elle-même: le
spectre de sa fausseté contre-exemplifiable fait partie de sa
rhétorique archaïque\,\,---\,\,inutile de rappeler cette
donnée\,\,---, et seuls les mathématiciens les plus avisés seront à
même de prendre à rebours les conjectures encore plus rares qui se
trompent d'orientation, parmi celles qui sont connues comme étant d'un
enjeu central\footnote{\, Le conjectural alimentaire de la recherche
mathématique courante est ici tenu à l'écart de
l'argumentation. }. Pour la même raison, il est fort inapproprié
d'écrire\,\,---\,\,même en acceptant l'intrusion de points
sophistiques involontaires\,\,---\,\,que la démonstration
<<\,\blue{\small\sf exclut que quelque chose puisse être appelé un
contre-exemple}\,>>: inapproprié en effet premièrement parce que les
démonstrations mathématiques n'orientent presque jamais\footnote{\, La
démonstration courante du théorème des quatre couleurs offre un
exemple exceptionnellement riche de stratégie d'élimination
systématique de contre-exemples potentiels.} leurs parcours en
excluant des contre-exemples potentiels: leur structure manifeste en
général un caractère argumentatif direct; mais ce n'est pas tout, cela
est inapproprié aussi, deuxièmement et même d'un point de vue
wittgensteinien <<\,puriste\,>>, puisque, si l'on admet comme le
soutient Wittgenstein que la structure logique intrinsèque de toute
démonstration doit s'identifier au seul sens que l'on peut conférer à
l'énoncé qu'elle démontre, alors toutes les fois qu'une démonstration
ne procède pas en éliminant tous les contre-exemples imaginables (ce
qui arrive la plupart du temps), il est \underline{faux} qu'une
démonstration <<\,\blue{\small\sf exclut que quelque chose puisse être
appelé un contre-exemple}\,>>. \`A tout le moins, une démonstration
sanctionnée doit exclure toute recherche de contre-exemple à l'énoncé
précis qu'elle démontre, sans pour autant empêcher de réfléchir à
l'optimalité des hypothèses en recherchant des contre-exemples à des
énoncés légèrement modifiés plus ambitieux. Tout énoncé est
accompagné d'un horizon coprésent de virtualités indécises concernant
les hypothèses qui le constituent.

\subsection*{ Reprise sur le théorème des nombres
premiers}
{\it Lucidité parfaite sur le fait que le sens de la proposition
s'identifie au contenu de sa démonstration effective; rigueur sur
l'étrangeté irréductible de nature entre l'inductif et le déductif}:
tel semble être l'apport majeur que Wittgenstein exprime de manière
récurrente comme s'il s'agissait de sa <<\,crispation spéculative\,>>
principale sur les mathématiques.

\smallskip

Sur le même exemple arithmétique continué, voici une confirmation
historique des écarts temporels importants qui peuvent séparer les
preuves des conjectures: après des travaux de Riemann, Bertrand,
Tchebychev, Mertens et d'autres, ce ne fut qu'un siècle après les
premières observations de Gauss, en 1896, que l'hypothèse quantitative
de répartition fut dûment et rigoureusement démontrée par Hadamard et
de la Vallée Poussin, en utilisant les méthodes transcendantes de la
théorie des fonctions d'une variable complexe: {\it Le nombre $\pi (
n)$ d'entiers positifs $p \in \{ 1, 2, 3, \dots, n\}$ qui sont
premiers tend vers l'infini de la même façon que $\frac{ n}{ \text{\rm
log}\, n}$.} Autrement dit:
\[
\lim_{n\to\infty}\,
\frac{\pi(n)}{\frac{n}{{\rm log}\,n}}
=
1, 
\]
ce que l'on écrit parfois $\pi ( n) \sim \frac{ n}{ {\rm log}\, n}$.

Ainsi, la loi expérimentale pressentie se révèle correcte. Seule une
démonstration est à même de spécifier comme vraie cette estimation
quantitative. En l'occurence ici, les premières démonstrations étaient
longues et délicates. Raison est donc donnée au philosophe analytique
wittgensteinien,
car tester ou découvrir expérimentalement cette loi en examinant une
liste de nombres premiers avec l'aide des tables de logarithmes
constitue une suite de gestes et d'actes de pensée (assez simples et
plutôt répétitifs) qui n'ont absolument rien à voir avec des arguments
de preuve délicats.

\smallskip

Concluons localement ces considérations: en revenant à l'extrait cité
{\it supra} page~\pageref{fosse}, il s'agit encore et toujours du même
<<\,fossé conceptuel\,>>, ou plus exactement d'une {\it distinction
fondamentale}\, entre:

\begin{itemize}

\smallskip\item[{\bf 1)}]
les \'enoncés mathématiques (notamment en théorie des nombres) qui
sont conjecturés grâce à des tests expérimentaux, à des calculs
numériques effectués automatiquement, à des listes exhaustives de
nombres, {\it etc.}, et:

\smallskip\item[{\bf 2)}]
les démonstrations mathématiques rigoureuses. 

\end{itemize}\smallskip

Approfondissons cela. En quoi et pourquoi est-il presque toujours
beaucoup plus facile de formuler des conjectures expérimentales que de
trouver des démonstrations? Cette question est subtile. Commençons par
une conjecture simple qu'aucune théorie n'accompage.

\subsection*{Conjecture de Collatz} 
Considérons le procédé suivant, que l'on peut expérimenter sur de
nombreuses sites Internet.  \'Etant donné un nombre entier initial
arbitraire $n \geqslant 1$, le remplacer par $n/2$ s'il est pair, ou
par $3 n + 1$ s'il est impair; itérer ce calcul; observer que pour
toutes les valeurs de $n$ jusqu'à, disons 100 ($\sim 10^{ 18}$ en
2007), on redescend toujours à $1$ (suivi du cycle $1 \to 4 \to 2 \to
1$) après un certain nombre d'itérations. {\it Conjecturer qu'il en va
de même pour tout entier $n$}.

\smallskip

Il n'existe pas de <<\,recette mathématique\,>> plus simple. Par
exemple, pour $n = 6$, on obtient la suite $6$, $3$, $10$, $5$, $16$,
$8$, $4$, $2$, $1$.  On peut en remplir le ventre des
ordinateurs. Seule limite physique: la taille des données
stockées. Pour $n = 11$, quatorze itérations sont nécessaires; pour $n
= 27$, cent-onze (\!!), et les termes intermédiaires montent jusqu'à
4\,858, redescendent à 911, remontent à 9232, avant de redescendre à 1
en sursautant plusieurs fois. Où est la difficulté? Dans l'absence de
loi simple? Dans le chaotique?

En 1996, T.~Oliveira e Silva a écrit un programme en langage C qui
calcule les trajectoires de toutes les valeurs initiales $n$
inférieures à une limite donnée. Une fois lancé, le programme couvre
des intervalles de $2^{ 50}$ entiers. Sur un ordinateur d'une mémoire
vive de 266 MHz, en tenant compte de raffinements algorithmiques
suggérés par E.~Roosendaal (un concurrent international), près de 400
millions d'entiers (en moyenne) peuvent être traités à chaque
seconde. Ce test fut stoppé quand $100 \cdot 2^{ 50}$ fut atteint.

Depuis Juin 2004, les efforts de vérification ont repris. Les calculs
tournent depuis plus de trois ans. Il sont distribués sur une
vingtaine d'ordinateurs, utilisent des algorithmes révisés qui sont
trois fois plus rapides que les précédents, et permettent de couvrir
des intervalles de $2^{ 58}$ entiers (facteur d'amélioration: $2^8 =
256$). Au printemps 2007, Collatz est confirmé pour tous les entiers
jusqu'à $14\cdot 2^{ 58} \simeq 4\cdot 10^{ 18}$. Par ailleurs,
Collatz se vérifie rapidement\footnote{\, {\tt
did.math.uni-bayreuth.de/personen/wassermann/fun/3np1.html} } pour des
entiers $n \leqslant 10^{ 500}$ tapés au hasard sur un clavier, parce
que les calculs sont triviaux pour la machine.

\subsection*{ La maxime capitale de l'induction} 
Soit une conjecture {\it ouverte}\, quelconque ${\sf Cjct} (n)$
portant sur une quantité qui dépend d'un entier $n$ arbitraire. Voici
ce que rappelle la <<\,\green{\sl maxime capitale de l'induction}\,>>:
quelle que soit la hauteur impressionnante\,\,---\,\,$n \leqslant
3\,000\,000$, $n \leqslant 10^{ 18 }$, $n \leqslant 10^{ 20}$, {\it
etc.}\,\,---\,\,jusqu'à laquelle ${\sf Cjct} (n)$ a été confirmée,
elle {\it peut toujours}\, être fausse. Sa probabilité de justesse,
comme sa probabilité de fausseté, sont essentiellement
inquantifiables.

\smallskip
\hfill
\begin{minipage}[t]{11cm}
\baselineskip=0.37cm
\blue{\small\sf
On a essayé d'évaluer la probabilité des inductions ou des hypothèses
en introduisant le concept de degré de confirmation d'une hypothèse
relativement à des faits. Ce degré de confirmation coïncide à peu
près avec une probabilité conditionnelle. Les logiques inductives que
l'on construit sur cette relation se sont révélées des formalismes
encombrants et inféconds. Il serait raisonnable de renoncer à trouver
à l'induction un fondement logique.
}
\hfill
{\small\green{Jean {\sc Largeault}.}}
\end{minipage}\medskip

L'indécision pure quant à la potentialité d'être ou de
ne pas être nous est imposée par l'imprévisibilité des mondes temporels. 
Misère et dénuement de l'entendement qui ignore!

Parfois, après des décennies de recherches, les réponses sont
crucifiantes. Plus d'une conjecture importante s'est révélée
contredite à des hauteurs exceptionnellement élevées de l'entier
$n$\footnote{\, Citons par exemple la conjecture de P\'olya, la
conjecture de Mertens et les nombres de Skewes. }.

\smallskip
\hfill
\begin{minipage}[t]{11cm}
\baselineskip=0.37cm
\blue{\small\sf
XXXV. \ \ \ \ \ \
Il n'est pas possible \`a celui qui commet clandestinement quelque
chose de ce que les hommes ont convenu entre eux de ne pas commettre
pour ne pas faire de tort ni en subir, d'\^etre s\^ur qu'il ne sera
pas d\'ecouvert, m\^eme si, dans le pr\'esent, il y \'echappe dix
mille fois, car, jusqu'\`a sa mort, l'incertain est s'il continuera
\`a n'\^etre pas d\'ecouvert.}
\hfill
{\small\green{{\sc \'Epicure},~{\em Maximes capitales}.}}
\end{minipage}\medskip

\medskip

Aussi l'évidence expérimentale ne {\it doit}-elle pas exister.
L'empiriste anti-inductif insiste: pour l'induction, il {\it doit}\,
ne pas y avoir de principe heuristique ou pseudo-probabiliste, {\it
parce que}\, le faux est toujours disponible dans l'ouvert. Le
philosophe analytique wittgensteinien navigue aussi dans ces prologues
de la spéculation mathématique spécialisée. 
Peut-il alors y avoir un
dogmatisme de l'indécision? 
\`A tout le moins, le {\it maintien
rigoureux de l'ouverture}\, constitue un {\it impératif catégorique}\,
de la pensée mathématique.

Mais la croyance en la véracité ou en la fausseté de ${\sf Cjct} ( n
)$ $\forall \, n$ doit forcer à engager des actes
irréversibles. Encore une fois: s'orienter, se confronter, c'est se
potentialiser, donc s'imprévisibiliser. Il va ainsi dans le monde
mathématique. 

Pour la conjecture de Collatz (ouverte depuis 1937), aucun appareil
théorique n'existe: c'est un cas exceptionnel. La sonde innocente:
\[
n \longmapsto \left\{
\aligned
&
n/2\ \
\text{\rm si}
\ \ n\ \ 
\text{\rm est pair}
\\
&
3n+1\ \
\text{\rm si}
\ \ n\ \ 
\text{\rm est {\it im}pair}
\endaligned
\right.
\]
est lancée dans l'indéfini potentiel primordial $\N = \{ 1, 2, 3,
\dots, n, \dots \}$. Il y a une {\it règle de calcul}, au sens de
Wittgenstein. Mais aucun encadrement théorique n'est connu, y 
compris pour d'autres sondes analogues\footnote{\, 
{\sc Lagarias}, J.C.: {\em The $3 x + 1$ problem
and its generalizations}, Amer. Math. Monthly
{\bf 92} (1985), 3--23.
}.

\smallskip
\hfill
\begin{minipage}[t]{11cm}
\baselineskip=0.37cm
\blue{\small\sf
Mathematics is not yet ready for such problems.
}
\hfill
{\small\green{Paul {\sc Erd\"os}.}}
\end{minipage}\medskip

Quel contraste entre cette règle d'itération simplissime et le chaos
des résultats obtenus! L'<<\,écart\,>>, le <<\,fossé conceptuel\,>>
se fait d'autant plus sentir qu'aucune démonstration n'existe en
germe. On ne dispose que d'un raisonnement probabiliste non rigoureux
pour se convaincre d'une éventuelle véracité de cette
conjecture\footnote{\, Si l'on considère seulement les nombres impairs
dans la suite de Collatz, alors {\it en moyenne} le nombre impair
suivant est multiplié par $3/4$. Voici l'argument heuristique.

Prenons un entier $n_0$ {\it impair} et itérons le procédé de Collatz
jusqu'à obtenir un prochain entier impair $n_1$. Que vaut en moyenne
le rapport $n_1 \big/ n_0$? En supposant que le devenir est soumis à
des lois probabilistes équidistribuées et mélangeantes, on a: une fois
sur deux $n_1 = (3 n_0 + 1) \big/2$; une fois sur quatre $n_1 = (3 n_0
+
1)\big/ 2^2$; une fois sur huit $n_1 = (3 n_0 + 1) \big/ 2^3$; {\it
etc.}; par conséquent, la croissance moyenne de taille attendue entre deux
entiers impairs consécutifs $n_0$ et $n_1$ devrait
être égale à:
\[
\Big(
\frac{3}{2}
\Big)^{1/2}
\,
\Big(
\frac{3}{2^2}
\Big)^{1/4}
\,
\Big(
\frac{3}{2^3}
\Big)^{1/8}
\cdots
=
\frac{3}{4}
<
1.
\]
Ainsi, cet argument suggère qu'en moyenne, les itérés impairs
décroissent d'un facteur $\frac{ 3}{ 4}$. Mais l'hypothèse
d'équidistribution et de mélange n'a pas encore pu être démontrée; de
plus, comme le raisonnement est probabiliste, même s'il était
rigoureux, il ne pourrait pas exclure l'existence de cycles élevés qui
seraient exceptionnels par rapport au comportement moyen. }.

S'il doit y avoir des lois prédisant le comportement de ces suites, il
faut les extraire du chaos expérimental. Deux lois conjecturales ont
été observées\footnote{\, Sur la page: {\tt
www.ieeta.pt/$\sim$tos/3x+1.html}, le lecteur trouvera deux graphiques
convaincants. }. Elles raffinent la perception de ce problème, sans
donneraucune indication de preuve.

\begin{itemize}

\smallskip\item[$\bullet$]
L'\green{\sl excursion maximale}\, de $n$, à savoir la valeur entière
maximale de sa trajectoire (9\,232 pour $n = 27$) semble se comporter
asymptotiquement comme $n^2$, tout en fluctuant autour de cette
valeur.

\smallskip\item[$\bullet$]
Le \green{\sl temps d'arrêt}\, d'un entier $n$, à savoir le plus petit
nombre d'itérations nécessaires pour passer en-dessous de $n$ (et se
ramener, par récurrence à un entier déjà examiné), semble se comporter
comme ${\rm log} \, n$, avec des fluctuations plus importantes.

\end{itemize}\smallskip

\subsection*{ Libération par le contre-exemple?}
Phénomène radicalement irréversible, l'avènement d'un contre-exemple
libère immédiatement de la question initialement posée, il libère d'un
travail de calcul indéfini, il arrête net une poursuite aveugle du
programme. \`A cet instant, toutes les intentions doivent changer,
tous les projets doivent être réorientés, toutes les intuitions être
réorganisées, et on stoppe les 20 ordinateurs calculant en parallèle
depuis plus de trois ans, et les 50 chercheurs concernés dans le monde
se remettent en question. C'est cela l'<<\,\green{\sl irréversible
mathématique}\,>>.

\subsection*{ Virtualités pérennes du principe de raison}
Mais très souvent, le contre-exemple révélé ne libère en rien de la
question en tant que question, parce que la question ne s'était
qu'imparfaitement exprimée dans la conjecture. La conjecture
prétendait que les êtres qu'elle interrogeait jouissaient d'une
certaine simplicité comportementale encadrée par certaines loi
quantitatives, mais elle n'effaçait pas toutes les complexités
adventices de ces êtres qui s'étaient déjà pré-exprimées dans les
moments de virtualisation collatérale. 

Le conjectural commence toujours par prétendre pour lui-même que le
simple domine, en tant que forme d'ensemble des phénomènes. Puis,
s'il se trompe, il corrige, il affine, il repousse, il accepte, il
complexifie. Curieusement, la dialectique du conjectural ne cesse de
remobiliser le même mouvement inépuisable de pensée qui cherche à
prévoir et à deviner des lois mathématiques régulatrices. Si le
<<\,principe de raison\,>> a envahi la pensée technicienne, comme l'a
parfois déploré Heidegger, cela même reste un mystère pour nous de
savoir ce qu'il y reste de pensée métaphysique et comment cette pensée
métaphysique irrigue encore continûment la pensée technique.

\subsection*{ Actifier la question}
Supposons découvert un cycle de Collatz très élevé\,\,---\,\,un
contre-exemple\,\,---\,\,mais faisons rigoureusement abstraction des
questions nouvelles qui surgiraient après coup. Le chaos stochastique
des boucles attirées par le cycle $1 \to 4 \to 2 \to 1$ {\it que l'on
avait déjà observé avant l'avènement dudit contre-exemple}\, n'en
serait pas moins mystérieux, {\it toujours}\, en question. Les
questions demeurent parce que le questionnement focalise son faisceau
sur des affirmations hypothétiques transitoires. Mais le
questionnement est toujours déjà éclaté au moment où il s'exprime. Le
questionnement est un acte spécifique de décision multiple que l'on
peut toujours reproduire, exporter, ramifier et faire éclater. {\it Le
questionnement est un \underline{acte} élémentaire}, un acte naturel
qui va de soi, et cet acte est analogue dans sa finitude aux actes
argumentatifs élémentaires du discours déductif (nous y
reviendrons). La forme même des questions mathématiques est
essentiellement simple, atomique.

\medskip

\hspace{-3cm}
\fbox{\red{\small\sf
Pourquoi la question mathématique atomique
produit-elle de l'irréversible-synthétique organique?}}

\medskip

Ici transparaît une thèse de philosophie des mathématiques que nous
jugeons capitale mais que nous ne dévoilerons pas encore pour
l'instant.

\smallskip
\hfill
\fbox{
\begin{minipage}[t]{11cm}
\baselineskip=0.37cm
\blue{\small\sf
We have to assume that we are very stupid and our natural questions
are stupid, and only by hard work, by conceptualizing, working hard,
calculating, whatever, we can make good questions or good mathematics.
And it's naive to think that we all have intuition or something. It's
a stupid opinion. That's what I believe.}

\hfill
{\small\green{Mikhail {\sc Gromov}.}}
\end{minipage}}\medskip

\subsection*{Conjecture de Proth-Gilbreath.} 
Deuxième exemple de conjecture purement expérimentale sans arrière
plan théorique dont on peut abreuver les ordinateurs. Voici la
recette, attribuée à Gilbreath, mais qui remonte à Proth au $\text{\rm
19}^{\text{\rm ième}}$ siècle. L'idée consiste à <<\,dévisser\,>> la
complexité des nombres premiers en calculant leurs différences
successives, ce à quoi Euler s'était adonné avec succès sur 
de multiples exemples.

\'Ecrire les nombres premiers les uns à la suite des autres sur une
première ligne; écrire sur une seconde ligne la valeur absolue des
différences entre deux nombres consécutifs; itéter cette opération;
conjecturer que chaque ligne commence par $\red{\bf 1}$:

\medskip

\begin{footnotesize}

\hspace{-2cm}
\begin{tabular}[t]{*{18}c}
\green{${}_k\diagdown {}^n$}&
&\green{1}&\green{2}&\green{3}&\green{4}
&\green{5}&\green{6}&\green{7}&\green{8}
&\green{9}&\green{10}&\green{11}&\green{12}
&\green{13}&\green{14}&\green{15}&\green{16}
\\
\cline{3-18}
\\
\green{0}&\vline
&{\bf 2}&{\bf 3}&{\bf 5}&{\bf 7}
&{\bf 11}&{\bf 13}&{\bf 17}&{\bf 19}
&{\bf 23}&{\bf 29}&{\bf 31}&{\bf 37}
&{\bf 41}&{\bf 43}&{\bf 47}&{\bf 53}
\\
\hspace{30mm}
\\
\green{1}&\vline
&\red{\bf 1}&\blue{2}&\blue{2}&\blue{4}
&\blue{2}&\blue{4}&\blue{2}&\blue{4}
&\blue{6}&\blue{2}&\blue{6}&\blue{4}
&\blue{2}&\blue{4}&\blue{6}
\\
\green{2}&\vline
&\red{\bf 1}&\blue{0}&\blue{2}&\blue{2}
&\blue{2}&\blue{2}&\blue{2}&\blue{2}
&\blue{4}&\blue{4}&\blue{2}&\blue{2}
&\blue{2}&\blue{2}
\\
\green{3}&\vline
&\red{\bf 1}&\blue{2}&\blue{0}&\blue{0}
&\blue{0}&\blue{0}&\blue{0}&\blue{2}
&\blue{0}&\blue{2}&\blue{0}&\blue{0}
&\blue{0}
\\
\green{4}&\vline
&\red{\bf 1}&\blue{2}&\blue{0}&\blue{0}
&\blue{0}&\blue{0}&\blue{2}&\blue{2}
&\blue{2}&\blue{2}&\blue{0}&\blue{0}
\\
\green{5}&\vline
&\red{\bf 1}&\blue{2}&\blue{0}&\blue{0}
&\blue{0}&\blue{2}&\blue{0}&\blue{0}
&\blue{0}&\blue{2}&\blue{0}
\\
\green{6}&\vline
&\red{\bf 1}&\blue{2}&\blue{0}&\blue{0}
&\blue{2}&\blue{2}&\blue{0}&\blue{0}
&\blue{2}&\blue{2}
\\
\green{7}&\vline
&\red{\bf 1}&\blue{2}&\blue{0}&\blue{2}
&\blue{0}&\blue{2}&\blue{0}&\blue{2}
&\blue{0}
\\
\green{8}&\vline
&\red{\bf 1}&\blue{2}&\blue{2}&\blue{2}
&\blue{2}&\blue{2}&\blue{2}&\blue{2}
\\
\green{9}&\vline
&\red{\bf 1}&\blue{0}&\blue{0}&\blue{0}
&\blue{0}&\blue{0}&\blue{0}
\\
\green{10}&\vline
&\red{\bf 1}&\blue{0}&\blue{0}&\blue{0}
&\blue{0}&\blue{0}
\\
\green{11}&\vline
&\red{\bf 1}&\blue{0}&\blue{0}&\blue{0}
&\blue{0}
\\
\green{12}&\vline
&\red{\bf 1}&\blue{0}&\blue{0}&\blue{0}
\\
\green{13}&\vline
&\red{\bf 1}&\blue{0}&\blue{0}
\\
\green{14}&\vline
&\red{\bf 1}&\blue{0}
\\
\green{15}&\vline
&\red{\bf 1}
\end{tabular}
\end{footnotesize}

\medskip
Plus précisément, soit $p_1 = 2$, $p_2 = 3$, $p_3 = 5$, $p_4 = 7$,
$p_5 = 11$, \ldots, les nombres premiers listés dans l'ordre croissant
et posons:
\[
\left[
\aligned
d_0(n)
&
:=
p_n,
\ \ \ \ \ \ \ \ \ \ \ \ \ \ \ \ \ \ \ \ \ \ \ \ \ \ \ \ \ \ \ \ \ \ \
n\geqslant 1,
\\
d_{k+1}(n)
&
:=
\big\vert
d_k(n)-d_k(n+1)
\big\vert,
\ \ \ \ \ \ \
k\geqslant 0,\ \ \ 
n\geqslant 1.
\endaligned\right.
\]
Le tableau montre que $d_k ( 1) = 1$ pour $1 \leqslant k \leqslant
15$. En 1959, la conjecture $d_k ( 1) = 1$ pour tout $k$ a été
confirmée par Killgrove et Ralston jusqu'aux profondeurs $k\leqslant
63\,419$, et pour tous les entiers premiers $< 792 \, 731$. En 1993,
A.M.~Od\-lyzko\footnote{\, {\it Iterated absolute values of
differences of consecutive primes}, Math. Comp. {\bf 61} (1993),
no.~203, 373--380. } confirme le phénomène pour tous les entiers
premiers $< 10^{ 13}$, de telle sorte que $d_k ( 1) = 1$ jusqu'à une
profondeur $\lesssim 3,4 \cdot 10^{ 11}$.

Pour une première ligne $d_0 ( n)$ qui serait
constituée d'entiers quelconques,
le calcul de $d_k ( 1)$ requiert en général {\it a priori}\, la
connaissance de tous les $d_j ( i)$ pour $i + j \leqslant k + 1$, de
sorte que pour $k \sim 3, 4 \cdot 10^{ 11}$, il faudrait calculer
aproximativement $5 \cdot 10^{ 22}$ nombres\,\,---\,\,au-delà des
capactités technologiques actuelles. Mais pour une première ligne
constituée des nombres premiers $d_0 ( n ) = p_n$, le tableau montre
qu'après un temps assez court, il n'y a plus que des 0 et des 2 après
le \red{\bf 1} attendu, et dans un tel cas, {\it i.e.} si, pour un $N$
on peut trouver un $K$ avec $d_1 ( 1) = \cdots = d_K (1) = 1$
tel que $d_K ( n) = 0$ ou $2$ pour tout $2
\leqslant n \leqslant N$, alors il est immédiat que l'on a ensuite
$d_k ( 1) = 1$ pour tout $K \leqslant k \leqslant N + K - 1$. Ce
phénomène se confirme et permet de réduire considérablement les temps
de calcul.

\smallskip
\hfill
\begin{minipage}[t]{11cm}
\baselineskip=0.37cm
\blue{\small\sf
A rigorous proof of Gilbreath's conjecture appears out of reach, given
our knowledge of primes. [\dots] About half of the machine time was
spent in sieving for primes, and half in computing the iterated
absolute values of the differences.
}
\hfill
{\small\green{A.M.~{\sc Odlyzko}.}}
\end{minipage}

\subsection*{ Retour sur Wittgenstein: deux universalités incomparables}
Reprenons l'opposition wittgensteinienne. 

\smallskip
\hfill
\begin{minipage}[t]{11cm}
\baselineskip=0.37cm
{\small{\sf\blue{
Ce qui n'est pas concevable, aux yeux de Wittgenstein, est que
l'universalité qui nous est fournie par la démonstration, lorsque nous
avons réussi effectivement à démontrer la proposition, puisse être
celle-là même que des expériences répétées, effectuées avec la méthode
de contrôle, nous avaient permis de supputer: <<\,Où est censée
ressortir de la démonstration la même universalité que les essais
antérieurs rendaient probables?>> (PG, p.~361.) Je peux assurément
formuler l'hypothèse douée de sens que, si je teste l'un après l'autre
les nombres pairs pour voir s'ils satisfont ou non la proposition de
Goldbach, je ne rencontrerai aucun contre-exemple de mon vivant. Mais
comment une démonstration de la proposition dans laquelle il n'est
question ni de moi, ni de qui que ce soit, ni de ce que je ferai ou ne
ferai pas, pourrait-elle démontrer cette supposition?
}}
\hfill
\green{\cite{ bo1988},~p.~194.}
}
\end{minipage}\medskip

Certainement, la démonstration ne ressort pas d'une série de tests
numériques. La nécessité universelle argumentée transcende la
confirmation expérimentale. Mais ici, encore une fois, on projette le
débat sur une opposition dualiste simplifiée. Alors que
l'irréversibilité historique de la mathématique impose une complexité
toujours grandissante aux dialectiques de la découverte, les
oppositions en restent ici à un stade non ramifié. L'histoire des
confirmations expérimentales s'étend sur plusieurs siècles; les
pratiques ont évolué; et l'ontologie physique du calcul s'est
considérablement enrichie à cause de la reproduction planétaire
des machines électroniques.  Atteindre un record de confirmation
expérimentale pour la conjecture de Goldbach n'a vraiment rien de
trivial actuellement; nous en reparlerons dans un instant.

De plus, l'affirmation <<\,\blue{\small \sf si je teste l'un après
l'autre les nombres pairs pour voir s'ils satisfont ou non la
proposition de Goldbach, je ne rencontrerai aucun contre-exemple de
mon vivant}\,>> part d'une prémisse insensée: aucun individu n'a jamais
consacré, et ne consacrera jamais l'intégralité de la durée de sa vie à
énumérer les cas d'une conjecture telle que celle de Goldbach les uns
à la suite des autres jusqu'à son dernier souffle.  Résumer sa vie à
une finitude éprouvée sur le parcours répétitif d'une seule
conjecture, ce serait se condamner et se crucifier. Mais en vérité,
nul ne songe à se priver du jeu de l'imprévu et du plaisir de décider
de ses propres changements d'orientation intellectuelle.

Autre objection: aujourd'hui, les confirmations expérimentales ne
s'effectuent plus à la première personne. Le <<\,je\,>> de l'activité
mathématique singulière n'a plus aucun sens, parce qu'il y a un
<<\,nous\,>> commun et international de la confirmation expérimentale,
qui tend de plus à se dépersonnaliser à cause de l'électronisation du
calcul, et de sa transmissibilité par les canaux de communication. Les
travaux de confirmation expérimentale se partagent entre les
chercheurs.

Poursuivons la critique.  Afin de défendre strictement sa thèse
dualiste du <<\,fossé\,>> entre expériences numériques et grammaires
formelles, le philosophe analytique wittgensteinien affirme que la
démonstration rigoureuse d'une proposition mathématique ne peut avoir
aucune incidence sur les suppositions qui se formulent {\it en tant
que telles}\, dans le champ de l'expérience.  Ou tout du moins, il
affirme que l'universalité hypothétique qui est suggérée dans
l'expérience n'est pas subsumée par l'universalité logique de la
démonstration, et partant, que l'universalité présumée conserve son
autonomie et son irréductibilité de principe. Cette affirmation
unilatérale est \underline{erronée}, et ce, pour quatre raisons.

\begin{itemize}

\smallskip\item[$\bullet$]
Parce qu'elle change le statut de la supposition expérimentale en
certitude universelle, la démonstration a un impact immédiat: le
caractère hypothétique, problématique et ouvert de la confirmation
disparaît du même coup, et toutes les tâches de vérification
calculatoire se métamorphosent en simples exercices d'application
numérique. 

\smallskip\item[$\bullet$]
Pour ce qui concerne la production et l'assimilation de
l'irréversible-synthétique, le principe de libre circulation entre
l'{\it a priori}\, et l'{\it a posteriori}\, exige que l'étudiant ou
le chercheur doive toujours pouvoir {\it se réinscrire temporairement
dans une situation d'ignorance artificialisée}\footnote{\, C'est bien
parce qu'on y efface les marques de l'indécision dialectique
originelle quant à l'irréversible-synthétique que les textes
mathématiques sont si difficiles à lire. }.

\smallskip\item[$\bullet$]
Dès qu'une conjecture est confirmée par une preuve, d'autres
suppositions plus ambitieuses peuvent être formulées en
partant de raisonnements heuristiques analogues.
L'homologie de structure se reproduit. 

\smallskip\item[$\bullet$]
La métaphysique audacieuse de la recherche entrelace tous
les niveaux formels et informels, avec toujours la même confiance
affirmée qu'il doit exister des lois et des démonstrations
potentielles.

\end{itemize}\smallskip

Certes, la nécessité apodictique de la démonstration ne {\it
démontre}\, pas quelque chose à propos des suppositions que nous
formulons en interrogeant les structures arithmétiques, ni même au
sujet de la mystérieuse faculté que nous avons d'énoncer de telles
suppositions, mais à tout le moins, il y a là un grand problème de
métaphysique des mathématiques qu'on ne peut pas se contenter
d'écarter obsessionnellement comme l'a fait Wittgenstein. L'optimisme de
Hilbert ({\it non ignorabimus}) et la méditation rétrograde de
Heidegger (domination universelle du principe de raison)
ressurgissent comme questions ouvertes de philosophie des
mathématiques.

\subsection*{ Exemple} 
Ainsi, nous affirmons que même après qu'une démonstration rigoureuse a
été produite, on peut exiger un retour vers l'expérimental numérique,
soit comme confirmation d'une sorte d'harmonie préétablie, soit comme
pénétration indépendante dans la réalité problématique
des mathématiques. Par exemple,
dans les années 1910 à 1920, G.~Hardy et S.Ramanujan\footnote{\, {\it
cf.} G.H.~{\sc Hardy}, {\em Some famous problems of the theory of
numbers and in particular Waring's problem. An inaugural lecture
delivered before the University of Oxford}, Oxford, Clarendon Press,
1920.}  ont découvert une formule approchée pour le nombre $p(n)$ de
partitions d'un entier $n$, dont le terme principal est:
\[
\frac{1}{2\pi\,\sqrt{2}}\,
\frac{d}{dn}\,
\frac{e^{\frac{2\pi}{\sqrt{6}}\,\sqrt{n- \frac{1}{24}}}}
{\sqrt{n-\frac{1}{24}}}.
\]

\smallskip
\hfill
\begin{minipage}[t]{11cm}
\baselineskip=0.37cm
\blue{\small\sf
[This formula] enables us to approximate to $p(n)$ with an accuracy
which is almost uncanny. We are able, for example, by using 8 terms of
our formula, to calculate $p ( 200)$, a number of 13 figures, with an
error of $0,004$. I have set out the details of the calculation:}
\end{minipage}\medskip
\[
\aligned
3\,972\,998\,993\,185,896
\\
36\,282,978
\\
-87,555
\\
5,147
\\
1,424
\\
0,071
\\
0,000
\\
0,043
\\
\hline
\\
3\,972\,999\,029\,388,004
\endaligned
\]
\begin{minipage}[t]{11cm}
\baselineskip=0.37cm
\blue{\small\sf
The value of $p(200)$ was subsequently verified by Major MacMahon, by
a direct computation which occupied over a month.
}
\hfill
{\small\green{G.H.~{\sc Hardy}.}}
\end{minipage}\medskip

\subsection*{Conjecture de Goldbach}
Venons en maintenant à un autre exemple
célèbre de conjecture ouverte: 
{\it tout nombre entier pair $\geqslant 4$ est somme de deux nombres
premiers}. Plus précisément, pour tout entier pair $2n \geqslant 4$,
il existe $p$ et $q$ appartenant à l'ensemble $\mathcal{ P}$ des
nombres premiers tels que $2 n = p + q$.

\smallskip

Le principe initial de la confirmation expérimentale est extrêmement
simple: il suffit en principe de se constituer au préalable une liste
de tous les nombres premiers (en utilisant par exemple le crible
d'\'Eratosthène, ce qui expose à la question d'efficacité et aux
difficultés d'implémentation) jusqu'à une certaine grandeur, de les
additionner deux à deux et d'examiner si tous les nombres entiers $n$
sont ainsi obtenus jusqu'à une certaine grandeur.

Pour confirmer cela dans un intervalle d'entiers $[ a, b]$, deux
méthodes ont été utilisées. On doit trouver deux ensembles de nombres
premiers $\mathcal{ P}_1 \subset \mathcal{ P}$ et $\mathcal{ P}_2
\subset \mathcal{ P}$ tels que
\[
\big\{
2n:\,\,
a\leqslant 2n\leqslant b
\big\}
\subset
\mathcal{P}_1+\mathcal{P}_2
=
\big\{
p_1+p_2:\,\,
p_1\in\mathcal{P}_1,\,\,
p_2\in\mathcal{P}_2
\big\}.
\]
Fixons un entier $\delta \geqslant 1$. 
Dans la première méthode on choisit:
\[
\mathcal{P}_1
=
\big\{
p_1\in\mathcal{P}:\,\,
2\leqslant p_1\leqslant b-a+\delta
\big\},
\ \ \ \ \ \ \ \ \ \
\mathcal{P}_2
=
\big\{
p_2\in\mathcal{P}:\,\,
a-\delta\leqslant p_2\leqslant a
\big\}.
\]
Dans la seconde méthode on choisit:
\[
\mathcal{P}_1
=
\big\{
p_1\in\mathcal{P}:\,\,
2\leqslant p_1\leqslant \delta 
\big\},
\ \ \ \ \ \ \ \ \ \
\mathcal{P}_2
=
\big\{
p_2\in\mathcal{P}:\,\,
a-\delta\leqslant p_2\leqslant b
\big\}.
\]
Les calculs montrent que $\delta$ peut en fait être choisi très petit
par rapport à $b$ pour trouver au moins un couple $( p_1 , p_2 ) \in
\mathcal{ P}_1 \times \mathcal{ P}_2$ tel que $2n = p_1 +
p_2$. L'évidence numérique supportant cette conjecture est très forte,
car le nombre \green{$g ( 2n)$} de \green{\it partitions de Goldbach},
{\it i.e.}  de manière d'écrire $2 n = p + q$ avec $p, q \in \mathcal{
P}$ et $p \leqslant q$, croît rapidement avec $2n$.

La première méthode a été implémentée sur des ordinateurs dès les
années 1960. Parce qu'elle exige d'effectuer des tests de primalités
sur de grands intervalles d'entiers $[ a - \delta , b ]$, la seconde est
moins économique, mais elle seule permet d'accéder à la \green{\sl partition
de Goldbach minimale}\, d'un entier pair $2n$ quelconque, c'est-à-dire
au couple d'entiers premiers $\big(p_{\rm min} (2n), q_{\rm min}
(2n) \big)$ avec $p_{\rm min} (2n) \leqslant q_{\rm min} (2n)$
tels que {\it pour tout autre}\, partition de Goldbach $2 n = p + q$
avec $p \leqslant q$, on a $p_{ \rm min} < p$.

En vérité, la recherche des partitions de Goldbach minimales expose à
une difficulté imprévisible: lorsque $2n$ augmente régulièrement, les
entiers premiers $p_{\rm min}$ sautent de manière assez
chaotique\footnote{\, Le lecteur trouvera une représentation graphique
de la \green{\sl voile de Goldbach}\, à la page: {\tt
wardley.org/images/misc/goldbach/}}.  Par exemple, juste avant $100 =
\red{\bf 3} + 97$, on a $98 = \red{\bf 19} + 79$. Ce phénomène
pourra-t-il être embrassé dans une démonstration d'une longueur
raisonnable?  Sinon, pourra-t-on le contourner grâce à une structure
globale de l'ensemble des partitions de Golbach de tous les entiers
pairs?

\medskip

\begin{center}
\fbox{\red{\small\sf L'expérimental numérique force à éclater 
les questions.}}
\end{center}

\subsection*{ Calculs au front}
Depuis février 2005, T.~Oliveira e Silva (en compétition avec d'autres
concurrents internationaux) pilote une cinquantaine d'ordinateurs qui
travaillent en parallèle pour chercher la partition de Goldbach
minimale d'entiers pairs appartenant à des intervalles de longueur
$10^{ 12}$.  En avril 2007, $10^{ 18}$ a été atteint.  Les calculs
mémorisent le nombre de fois que chaque (relativement petit) nombre
premier $p$ est utilisé dans une partitition de Goldbach minimale,
ainsi que le plus petit entier pair $2n$ pour lequel $p_{\rm min} (2n)
= p$.

\smallskip
\hfill
\begin{minipage}[t]{11cm}
\baselineskip=0.37cm
\blue{\small\sf
On a 2.2GHz Athlon64 3500+ processor, testing an interval of $10^{
12}$ integers near $10^{ 18}$ takes close to 75 minutes. The execution
time of the program grows very slowly, like ${\rm log} (N)$, where $N$
is the last integer of the interval being tested, and it uses an
amount of memory that is roughly given by $13 \sqrt{ N} \big/ {\rm
log}\, N$. The program is now running on the spare time of around 50
computers (20 DETI/UA and 30 at PSU), either under GNU/Linux or under
Windows 2000/XP. We have reached $10^{ 18}$ in April 2007, and are now
double-checking a small part of the results.
}
\hfill
{\small\green{T.~{\sc Oliveira e Silva}}}
\end{minipage}\medskip

Pour tout entier premier $p$, les spécialistes se sont aussi
intéressés à \green{$S ( p) =$} le plus petit entier pair $2n$ tel que
$p$ apparaît dans la partition de Goldbach minimale de $2n$.  Le
record actuel (automne 2007) est détenu par J.~Fettig et N.~Sobh:
c'est $p = 9341$ pour $2n = 906\, 03057\, 95622\, 79642$.  En 1989,
A.~Granville, J.~van de Lune et H.J.J.~te Riele\footnote{\, {\em
Checking the Goldbach conjecture on a vector computer}, Number Theory
and Applications, R.A.~Mollin (ed.), pp.~423--433, Kluwer Academic
Press, 1989.} ont conjecturé, en invoquant un argument probabiliste,
que $p$ ne devait par croître plus rapidement que ${\rm log}^2 \, S (
p) \, {\rm log}\, {\rm log}\, S ( p)$.  Mais les données
expérimentales contredisent cette estimation qui devrait être
remplacée par $\frac{ 1}{ 3}\, \big[ {\rm log}\, S ( p) \, {\rm log}\,
{\rm log}\, S ( p) \big]^2$.

\medskip

\begin{center}
\fbox{\red{\small\sf L'expérimental numérique éprouve les cohérences
heuristiques.}}
\end{center}

\smallskip
\hfill
\begin{minipage}[t]{11cm}
\baselineskip=0.37cm
\blue{\small\sf
Explorer l'univers des nombres comme le monde physique?  Aucun domaine
n'a engendré autant de conjectures indémontrées (mais en partie
vérifiables à l'aide de calculs) que l'arithmétique des nombres
premiers. Contrairement à l'idée que les mathématiciens proposent le
plus souvent de leur discipline, les démonstrations y semblent parfois
reléguées au second plan.  De toute façon, disent les mathématiciens
eux-mêmes, nous n'arrivons pas à démontrer nos conjectures, et l'état
actuel de nos connaissances rend impensable que nous réussissions dans
un proche avenir.
}
\hfill
{\small\green{J.-P.~{\sc Delahaye}.}}
\end{minipage}\medskip

\subsection*{ Digression sur la nature physique du calcul}
Mais quelle {\it magie}\, alors nous délivrent les ordinateurs? Rien d'autre
qu'une mécanisation des gestes de calcul de type eulérien ou gaussien,
lorsque lesdits gestes s'astreignent {\it sans pensée latérale}\, à
aligner les résultats successifs obtenus par application d'une
certaine règle définie d'engendrement arithmétique.

\smallskip
\hfill
\begin{minipage}[t]{11cm}
\baselineskip=0.37cm
{\small{\sf\blue{
D'un bout à l'autre du calcul} [dans la preuve du théorème des quatre
couleurs]\blue{, n'importe qui peut étudier et vérifier chaque détail. Le
fait qu'un ordinateur puisse traiter en quelques heures plus de cas
particuliers qu'un humain ne pourrait espérer le faire dans toute sa
vie ne change rien au concept même de démonstration. 
}}
\hfill
\green{W.~{\sc Haken}.}
}
\end{minipage}\medskip

\noindent
L'ordinateur programmé par le théoricien expérimental des nombres
n'est donc rien de plus qu'un <<\,Train de calculs à Grande
Vitesse\,>> lancé dans l'indéfini primordial et irréductible qu'est la
suite des nombres entiers.

\smallskip
\hfill
\begin{minipage}[t]{11cm}
\baselineskip=0.37cm
{\small{\sf\blue{
Tous les calculs sont empiriques au sens trivial où ils supposent la
mise en {\oe}uvre d'une manipulation de symboles, que ce soit
mentalement, avec du papier et un crayon, ou à l'aide d'une machine!
}}
\hfill
\green{Martin~{\sc Gardner}.}
}
\end{minipage}\medskip

\noindent
Grâce aux microprocesseurs, la <<\,physicalité du calcul\,>> est ainsi
enrichie à un niveau micro- ou nano-scopique toujours plus
profondément lointain des bouliers orientaux, tables de calculs,
bâtonnets de Neper, machines à calcul mécaniques (Vinci, Schikard,
Pascal) ou machines à calcul électromécaniques, qui étaient
initialement conçues à l'échelle physique de l'homme. {\it En dernier
recours, les symboles en mouvement nécessitent toujours un support
matériel pour s'exécuter dynamiquement}. Les <<\,gestes de calcul\,>>
peuvent être compressés dans l'espace-temps et augmentés en volume:
telle est la seule et unique <<\,magie\,>> des ordinateurs. Et pour ce
qui est de l'essence même du calcul, la vraie et seule <<\,magie\,>>
qui nous entoure tous remonte aux babyloniens: c'est la
possibilité\,\,---\,\,au lieu de solliciter membres et
neurones\,\,---\,\,de piloter cailloux ou électrons dans l'univers
mobile du monde physique pour que ces éléments physiques calculent
automatiquement.

\medskip\noindent
{\sc Physicalité fondamentale du calcul.} {\it Qu'il soit manuel ou
digital, arithmétique, algébrique, numérique, probabiliste ou
diagrammatique, tout calcul exécuté ou programmé par les hommes est
irréductiblement discret, fini et imprimé de manière transitoire sur
des supports physiques. Aucun calcul <<\,transcendant\,>> à une
effectuation incarnée physiquement n'est possible. Tous les calculs
pour lesquels l'ordinateur est très performant (décimales de $\pi$;
bases de Gr\"obner; tests de primalité; analyse matricielle; schémas
numériques des équations aux dérivées partielles; statistiques; tris
de données) sont dans leur principe effectif identiques à
ceux que l'on conduit en ayant recours à n'importe quel autre
véhicule physique pour le mouvement des symboles.
}

\medskip

(Il reste toutefois très incertain que la puissance des ordinateurs
soit sans conteste effectivement supérieure à celle d'un Euler ou d'un
Gauss, même envisagés artificiellement comme n'étant que calculateurs
de génie: nous y reviendrons en temps voulu. Par ailleurs, il existe
de nombreux domaines des mathématiques qui ne se prêtent absolument
pas à une <<\,physicalisation\,>>, ni à aucun type
d'assistance électronique.)

\subsection*{ Dialectique a priori de l'existentiel}
L'atomicité symbolique du quantificateur <<\,$\exists$\,>> qui sert à
exprimer conjectures et théorèmes dans le même langage formel ne doit
pas faire croire que l'existence se réduise à un concept non
problématique. En mathématiques, l'existence ouverte est d'une
complexité dialectique imprévisible; ses variations spéculatives
peuvent s'avérer troublantes.

Rappelons que le débat philosophique entre l'existence abstraite et
l'existence effective en mathématiques (formalistes contre
constructivistes, Hilbert contre Gordan) est causé, en amont des
polémiques, par le fait que certains énoncés mathématiques peuvent
souvent être jugés comme imparfaits, partiels et donc encore ouverts
du point de vue de la connaissance mathématique\footnote{\, Nous
développerons cette thèse en temps voulu. }.

Ici\,\,---\,\,phénomène surprenant et paradoxal\,\,---, la conjecture
de Goldbach montre qu'{\it un trop-plein d'existence peut faire
obstacle à une connaissance mathématique achevée}: l'expérience montre
en effet que le nombre de couples de nombres premiers $(p_1, p_2)$
tels que $n = p_1 + p_2$ augmente très rapidement avec $n$.  La
dialectique {\it a priori}\, de l'existentiel ouvert doit donc
s'enrichir de ce cas de figure, et le considérer comme
métaphysiquement disponible à l'avenir.

\subsection*{ Heuristique semi-rigoureuse}
En 1923, grâce à des arguments informels mais pertinents, Hardy et
Littlewood ont conjecturé que le nombre $\pi_2 ( n)$ de
représentations de tout entier $n$ assez grand comme somme de deux
nombres premiers $n = p_1 + p_2$ devait être asymptotiquement égal
à:
\[
2\,\varpi_2\,
\frac{n}{\big({\rm log}\,n\big)^2}\,
\prod_{p\vert n;\,\,p\geqslant 3}\,
\frac{p-1}{p-2},
\]
où $n$ est pair et où $\varpi_2$ est la \green{\sl constante des
nombres premiers jumeaux}: $\varpi_2 = \prod_{ p \geqslant 3}\, \Big(
1 - \frac{ 1}{ ( p-1)^2} \Big) = 0,66016 \cdots$. Cette valeur
asymptotique est bien confirmée jusqu'à $n \leqslant 10^{ 17}$.  Les
manuels ou digitaux confirment la présence du facteur $\prod_{ p \vert
n;\,\,p\geqslant 3}\,\frac{p-1}{p-2}$ découvert par Sylvester en 1871
et qui produit de petites oscillations dans la valeur expérimentale de
$\pi_2 ( n)$ lorsque $n$ varie. Nous y reviendrons.

\smallskip

Considérons maintenant quelques conjectures ou 
questions ouvertes en arithmétique des nombres
premiers qui sont
simples à comprendre et à énoncer.

\subsection*{ Conjecture des nombres premiers jumeaux:} 
{\it Il existe un nombre infini de paires de nombres premiers $( p, p+
2)$ séparés seulement par un écart de $2$}. Autrement dit:
\[
\liminf_{n\to\infty}\,
p_{n+1}-p_n
=
2.
\]

\medskip

\subsection*{ Conjecture de Polignac} 
{\it Pour tout écart pair $2k$, il existe une
infinité de paires de nombres premiers $(p, p+ 2k)$ séparés par $2k$}.

\subsection*{ 
Existe-t-il une infinité de nombres premiers de la forme $n^2 +
1$?} On sait qu'il en existe une infinité de la forme $n^2 + m^2$ ou
$n^2 + m^2 + 1$.

\subsection*{ 
Existe-t-il toujours un nombre premier entre $n^2$ et $(n+
1)^2$?} En 1882, Opperman conjectura que $\pi ( n^2 + n) > \pi ( n^1)
> \pi ( n^2 - n)$, ce qui est aussi très probable.

\subsection*{ 
\'Ecarts entre nombres premiers consécutifs} En 1936, Cramér
conjectura que
\[
\limsup_{n\to\infty}\,
\frac{p_{n+1}-p_n}{\big({\rm log}\,p_n\big)^2}
=
1,
\]
d'où en particulier: il existe des écarts arbitrairement grands entre
nombres premiers qui se suivent.

\subsection*{ 
Fréquence des \'ecarts entre nombres premiers consécutifs} Wolf,
Odlyzko et Rubinstein ont conjecturé que l'écart le plus fréquent
entre deux nombres premiers est égal au produit des $n$ premiers
nombres premiers
\[
{\sf E}(n)
:=
2\times 3\times 5\times 7\times\cdots\times p_n
\]
pour tous les nombres compris entre
\[
{\sf h}(n)
:=
\exp
\Big(
\frac{2\times 3\times\cdots\times p_{n-1}(p_n-1)}{
{\rm log}\big[
(p_n-1)\big/(p_n-2)
\big]}
\Big)
\]
et ${\sf h} (n+1)$. Ici, ${\sf h}(3) \simeq 10^{ 36 }$ est déjà bien
au-delà du domaine accessible à une expérimentation
systématique\footnote{\, Contrairement aux conjectures précédentes,
l'expérimentation numérique ne peut pas constituer ici la source
principale d'alimentation prospective. Le dispositif initial du
philosophe analytique wittgensteinien est donc spéculativement
incomplet: il faut aussi tenir compte des conjectures qui hybrident un
champ expérimental insuffisant à des raisonnements heuristiques
semi-rigoureux. }, et {\it a fortiori}\, aussi ${\sf h} (4) \simeq
10^{ 428}$, ${\sf h} (5) \simeq 10^{ 8656}$, {\it etc}.

\subsection*{ Métaphysique du <<\,tout ce qui est possible se réalise\,>>}
Dans le domaine des nombres premiers, on peut formuler de très
nombreuses conjectures simples au sujet d'ensembles de nombres
astreints à satisfaire un certain nombre de propriétés définies. Tout
le possible qui n'est pas exclu par des conditions nécessaires
raisonnables et évidentes semble pouvoir prétendre à une plénitude
d'être, à une infinitude, à une quantifiabilité explicite. Face à la
réalité problématique irréductible des nombres entiers, et bien qu'il
semble ne pas y avoir de principe supérieur pour expliquer comment les
réalisations mathématiques sont possibles, la posture métaphysique du
conjectural engageme vers les potentialités du possible.  Parce que
l'expérience acquise par l'histoire des mathématiques témoigne de
réussites passées, les actes conjecturaux sont généralisables,
universalisables et reproductibles. Les formes du questionnement
mathématique s'organisent en une algèbre libre, ouverte et non
systématisable.

Mais d'un autre côté, la conjecture n'est qu'une forme d'accès
préliminaire aux réalités problématiques des mathématiques. Les formes
abstraites générales de l'interrogation exposent à de
l'irréversible-synthétique qui exige une circulation permanente 
des questions dans les preuves.

\smallskip
\hfill
\begin{minipage}[t]{11cm}
\baselineskip=0.37cm
\blue{\small\sf
The achievement of the mathematicians who found the Prime Number
Theorem was quite a small thing compared with that of those who found
the proof. [\dots] The whole history of the Prime Number Theorem, and
the other big theorems of the subject, shows that you cannot reach any
real understanding of the structure and meaning of the theory, or have
any sound instincts to guide you in further research, until you have
mastered the proofs. It is comparably easy to make clever guesses;
indeed there are theorems, like the ``Goldbach's theorem'', which have
never been proved and which any fool could haved guessed.
}
\hfill
{\small\green{{\sc G.~H.~Hardy.}}}
\end{minipage}

\subsection*{ Raisonnement absurde}
Reprenons maintenant notre analyse critique des expressions qui sont
employées par le philosophe analytique wittgensteinien.  Voici un
autre extrait.

\smallskip
\hfill
\begin{minipage}[t]{11cm}
\baselineskip=0.37cm
\blue{\small\sf
Cette idée qu'il existe une différence de nature, et non pas
simplement de degré, entre la démonstration et l'expérience, qui fait
que la démonstration ne peut pas démontrer exactement ce qui a été
conjecturé ({\it sic}), est liée au fait que, dans la proposition
mathématique, l'expression <<\,nécessairement tous\,>> constitue pour
ainsi dire un mot unique ({\it cf.} {\sc pg}, p.~429) et que l'on ne
peut en détacher le <<\,tous\,>> pour le comparer à celui de
l'expérience. Supposer que tous les nombres naturels ont une certaine
propriété veut dire supposer que, si on les passait tous en revue
successivement, on constaterait que chacun d'entre eux a cette
propriété. Mais que peut vouloir dire supposer que tous les nombres
naturels ont {\it nécessairement}\, une certaine propriété, si ce
n'est précisément supposer l'existence d'une démonstration de la
proposition universelle?
}
\hfill
{\small\green{\cite{bo1988},~p.~195.}}
\end{minipage}\medskip

Ici, la spéculation dérape: aveuglée par le même et unique dipôle
conjectures\big/preuves, elle exagère les différences conceptuelles en
extrapolant la signification de l'écart. Ici, {\it la tentation sophistique
menace l'exégète}. Même en admettant que les énoncés visés se
métamorphosent souvent au cours d'une recherche, et donc que les
démonstrations ne démontrent pas toujours nécessairement ce qui a été
initialement conjecturé ou visé, l'affirmation brutale <<\,\blue{\small\sf la
démonstration ne peut pas démontrer exactement ce qui a été
conjecturé}\,>> est \underline{inadmissible}:

\begin{itemize}

\smallskip\item[$\bullet$]
Sans autre nuance restrictive que par insertion furtive de l'adverbe
<<\,\blue{\small\sf exactement}\,>>, cette affirmation péremptoire se
présente comme valable pour {\it toute}\, proposition conjecturée et
pour et {\it toute}\, démonstration! Indignation chez les
mathématiciens!

\smallskip\item[$\bullet$]
Par une sorte d'argument d'autorité philosophique, cette affirmation
semble suggérer que celui qui démontre est toujours na\"{\i}f de
croire que ce qu'il démontre est effectivement ce qu'il annonce comme
ce qu'il va démontrer.

\smallskip\item[$\bullet$]
De plus, cette affirmation élimine brutalement tout ce qui fait
l'intention d'un projet déductif.

\smallskip\item[$\bullet$]
Enfin, plus grave encore, par l'insertion du verbe modal
<<\,\blue{\small\sf peut}\,>>, cette affirmation catégorique se
présente comme une vérité de fait qui limite {\it a priori}\, la
portée de toute démonstration par rapport à un énoncé.

\end{itemize}\smallskip

Et pour justifier cette absurde affirmation, en faisant un appel
rhétorique distendu et indirect à la périphrase imprécise
<<\,\blue{\small\sf est liée au fait que}\,>>, on greffe un appel au
quantificateur logique universel afin de convaincre définitivement son
lectorat de philosophie analytique: en tant qu'il est porteur
d'une nécessité logique, le quantificateur universel <<\,$\forall$\,>>
transcende le caractère inductif de la conjecture.

Ensuite, l'obsession portant sur le <<\,fossé conceptuel\,>> entre
expériences et preuves conduit à écrire une phrase surprenante:
<<\,\blue{\small\sf Supposer que tous les nombres naturels ont une
certaine propriété veut dire supposer que, si on les passait tous en
revue successivement, on constaterait que chacun d'entre eux a cette
propriété}\,>>: éh bien justement non! Sauf de manière accessoire et
partielle, ce n'est vraiment pas d'une vérification indéfinie que
parle une supposition mathématique! Et ce, pour deux raisons.

Premièrement, à l'échelle humaine, l'infini n'existe pas; on ne peut
jamais supposer qu'une infinité de nombres entiers soient passée en
revue: cela n'existe pas; cela ne peut pas exister. La thèse lucide
sur la physicalité du calcul que nous avons énoncée il y a quelques
instants a pour conséquence immédiate de borner (disons par $10^{
70}$) le nombres d'opérations jamais effectuables dans l'univers.

Deuxièmement, qu'elles soient conjecturales-ouvertes,
conjecturées-éta\-blies, ou simplement établies-admises, la plupart des
propositions mathématiques s'expriment dans un langage logique, avec
des quantificateurs existentiels ou universels. Il est lointain, le
temps où le langage axiomatique balbutiait! 

\smallskip

\begin{center}
\fbox{\red{\small\sf Le conjectural s'inscrit
d'emblée dans le langage du démonstratif.}}
\end{center}

\smallskip

En mathématiques,
l'universalité et l'existentialité du conjectural sont du même type 
métaphysique que
l'universalité et l'existentialité du démonstratif. La différence entre
les deux est toute modale: elle a trait au caractère d'{\it
ouverture}\, des énoncés. Bien que la tradition classique de
philosophie des mathématiques semble s'être résolument écartée de
l'ouverture comme concept, de grands mathématiciens comme Riemann ou
Hilbert nous ont légué quelques précieuses pensées à ce sujet. Nous y
reviendrons ultérieurement.

Poursuivons la critique. En mathématiques, un raisonnement est absurde
lorsqu'il est contradictoire. Jusqu'à nouvel ordre, le principe de
non-contradiction doit être rigoureusement respecté. En philosophie
spéculative, notamment dans la {\it Weltanschaaung}\, hégélienne, on
admet que ce principe puisse être remis en cause. Mais en philosophie
analytique, il est encore considéré à juste titre comme exigence
minimale. Or dans cet extrait, la cohérence locale des raisonnements
n'est pas respectée, parce que l'universalité et l'existentialité de
la proposition mathématique conjecturale s'expriment la plupart du
temps dans un langage formalisé qui attend une démonstration complète
exprimée dans le même langage et qui est accompagnée de démonstrations
partielles, d'idées initiales, d'arguments heuristiques.

\subsection*{ Maintien du fossé conceptuel} Encore
une citation témoignant de la circularité de la spéculation.
Le commentaire critique est laissé en exercice.

\smallskip
\hfill
\begin{minipage}[t]{11cm}
\baselineskip=0.37cm
\blue{\small\sf
La tentation à laquelle il faut résister, en l'occurence, est celle
qui consiste à considérer une série d'expériences de mesure
susceptibles de conduire à l'idée du théorème de Pythagore et la
démonstration du théorème comme deux symptômes différent du même état
de chose, le deuxième ayant simplement sur le premier l'avantage
d'être beaucoup plus sûr, et pour tout dire, infaillible.
Wittgenstein réagit à ce genre de suggestion en remarquant que:
<<\,rien n'est plus funeste pour la compréhension philosophique que la
conception de la démonstration et de l'expérience comme étant deux
méthodes de vérification différentes, donc tout de même
comparables\,>> ({\sc pg}, p.~361).
}
\hfill
{\small\green{\cite{ bo1988},~p.~195.}}
\end{minipage}\medskip

\subsection*{ Retour sur le théorème des nombres premiers; 
doxas anachroniques} La
régularité $\pi ( n) \sim \frac{ n}{ {\rm log}\, n}$ est surprenante:
elle aurait tout aussi bien pu se révéler fausse, si l'on s'en était
tenu à l'exercice spéculatif universel que nous offre la {\it doxa}\,
pure et {\it a priori}\, du conjectural. \`A partir du moment où la
suite des nombres premiers est considérée comme irréductible à toute
saisie formelle totalisante parce qu'indéfiniment riche et complexe,
comment cette suite pourrait-elle jouir de régularités aussi simples
que $\pi ( n) \sim \frac{ n}{ {\rm log}\, n}$? Et si l'on doit
admettre que de telles régularités simples existent effectivement,
comment se constituer une intuition fiable des structures plausibles?

Voilà encore un autre type de question qui demeure toujours en suspens et
toujours disponible: {\it comment un résultat établi s'insère-t-il dans
l'intuition provisoirement constituée qu'on a d'un champ
rationnel}? 

Bien que l'irréversible-synthétique engendre ses raisonnements
rigoureux, il est totalement faux que l'{\it a posteriori}\,
démonstratif efface l'ouverture fondamentale qui est inhérente à la
proposition non démontrée. L'imparfait demeure et l'ouverture
latérale non écrite reste coprésente. La faculté d'interrogation est
intacte: dès lors qu'on cherche à comprendre une démonstration en
profondeur, on doit métamorphoser, retourner et dés-{\it a
postérioriser} tous les raisonnements. On doit faire ressurgir les
questions décisives qui ont orienté l'irréversible-synthétique vers la
mise au point d'arguments spécifiques. La consignation des résultats
mathématiques dans un langage formel élague des dialectiques
qu'il faut reconstituer.

Le penseur wittgensteinien se trompe donc ici sur un point crucial: le
temps de la pensée circule dans tous les sens et voyage de manière
anachronique entre l'{\it a priori}\, et l'{\it a posteriori}, entre
la démonstration actuelle et sa saisie comme horizon, même si l'irréversible
biologique et la flèche du temps imposent que ces voyages s'effectuent
au détriment du vieillissement corps, même si les répétitions, les
hésitations, les reprises, les corrections se déploient linéairement
dans un temps biologique irréversible. On ne fait jamais réellement
abstraction du fait qu'un énoncé dit quelque chose que l'esprit
embrasse aisément en un instant, alors que l'étude des
démonstrations exige en général des heures de concentration
et de réflexion.

De plus, la démonstration ne supplante jamais définitivement son champ
expérimental originaire. Quiconque est intéressé par la répartition
des nombres premiers aura avantage à reprendre les tests de Gauss, et
il découvrira, comme Gauss, des oscillations locales presque
chaotiques dans cette répartition, oscillations que le théorème $\pi (
n) \sim \frac{ n}{ {\rm log}\, n}$ est visiblement incapable de
quantifier et au sujet desquelles il ne dit rien. En extrayant les
régularités essentielles, le conjectural se focalise sur les grandes
structures. Mais l'expérimental latéralise, complexifie et ramifie les
intuitions questionnantes. Les lois ne sont pas données d'emblée avec
les listes expérimentales: tout ce que l'on peut dire, c'est qu'elles
y transparaissent {\it peut-être}. L'expérience scientifique confirme
toujours la déraisonnable efficacité du principe de raison. Formuler
une loi conjecturale requiert toujours un acte synthétique de l'esprit.

\subsection*{ Permanence du provisoire et de la problématicité}
Aussi le schéma simplifié <<\,\blue{\small\sf conjectures versus
preuves}\,>> que le philosophe analytique retient d'une lecture de
Wittgenstein ne correspond-il vraiment pas à la complexité des
situations de recherche que provoque l'interrogation expérimentale,
toujours ouverte à des phénomènes subsidiaires. Dire que des questions
nouvelles renaissent une fois les résultats acquis serait encore
insuffisant, parce que: 

\medskip
\hspace{-1cm}
\fbox{\red{\small\sf Les questions intrinsèques perdurent au
sein des architectures achevées.}}

\medskip\noindent
Au sein même des démonstrations
purifiées, l'ouverture se maintient dans les questions qui sont déjà
tranchées.

Par ailleurs, et d'une manière générale, dans la pratique
mathématique, il y a un certain nombre de questions universelles
reproductibles. Ici par exemple, au sujet de la preuve de type
Hadamard et de la Vallée-Poussin, quelques questions à caractère
essentiellement universel peuvent être posées:

\begin{itemize}

\smallskip\item[$\bullet$]
comment les nombres premiers s'intègrent-ils dans l'analyse complexe?

\smallskip\item[$\bullet$]
quels sont les arguments décisifs? comment les distinguer
des arguments élémentaires? 

\smallskip\item[$\bullet$]
quelles sont les intuitions globales survolantes 
de la preuve? 

\smallskip\item[$\bullet$]
{\it la démonstration que je lis constitue-t-elle
la <<\,bonne\,>> démonstration?}

\end{itemize}\smallskip

Rien de plus permanent et de plus ineffaçable que les questions de
compréhension, notamment en mathématiques. C'est parce que le
conjectural contient des traces indélébiles de problématicité qu'il
est ineffaçable.

Considérons par exemple la quatrième question. \`A ce jour,
essentiellement deux démonstrations de l'équivalence $\pi ( n) \sim
\frac{ n}{ {\rm log}\, n}$ sont connues. La première, due à Hadamard
et de la Vallée Poussin, utilise la fonction $\zeta ( s)$ de Riemann,
la théorie des intégrales, les séries et produits infinis,
l'intégration dans le champ complexe, l'étude des valeurs au bord de
fonctions holomorphes, et des arguments de type taubérien: elle n'est
décidément pas <<\,élémentaire\,>>. Hadamard et de la Vallée Poussin
ont d'abord démontré que $\zeta ( s)$ ne s'annule pas dans le
demi-plan fermé $\{ {\rm Re}\, s \geqslant 1\}$, et ensuite établi des
estimées techniques de croissance $\zeta ( s)$ en $\infty$, afin
d'intégrer sur certains contours de Cauchy allant à l'infini pour
obtenir les coefficients de séries de Dirichlet (comme $\zeta (
s)$). L'étude préliminaire de $\zeta ( s)$ a été simplifiée par
Tchebychev, Titchmarsh et Mertens. Le recours aux séries de Fourier
(Wiener, Ikehara, Heins) offre une alternatives aux arguments finaux
de Hadamard et de la Vallée Poussin. Mais actuellement, la preuve la
plus concise et la plus directe, qui n'utilise presque rien de plus
que la formule de Cauchy, a été mise au point par
D.~J.~Newman\footnote{\, Dans un article dédié au centième
anniversaire du théorème des nombres premiers qui est paru en 1997 à
l'{\it American Mathematical Monthly}, vol.~{\bf 10}, 705--708, Don
{\sc Zagier} restitue la preuve de Newman en trois pages d'une
limpidité et d'une concision remarquables. La preuve procède en six
moments. Pour $s\in \C$ et $x\in \R$, définissons:
\[
\zeta(s)
:=
\sum_{n=1}^\infty\,\frac{1}{n^s},
\ \ \ \ \ \ \ \ \ \ 
\Phi(s)
:=
\sum_p\,\frac{{\rm log}\, p}{p^s},
\ \ \ \ \ \ \ \ \ \ 
\vartheta(x):=\sum_{p\leqslant x}\,
{\rm log}\, p,
\]
où la lettre $p$ est utilisée pour désigner les nombres premiers.

\noindent
{\bf I:} $\zeta ( s) = \prod_p\, (1 - p^{ -s} )^{ -1}$ pour
${\rm Re}\, s > 1$.

\noindent
{\bf II:}
$\zeta (s) - \frac{ 1}{ s - 1}$ se prolonge
holomorphiquement à $\{ {\rm Re}\, s > 0 \}$.

\noindent
{\bf III:}
$\vartheta ( x) = {\rm O} ( x)$.

\noindent
{\bf IV:} $\Phi (s) - 
\frac{ 1}{ s - 1}$ est holomorphe et $\zeta (s)$ ne s'annule pas
dans
$\{ {\rm Re}\, s \geqslant 1\}$.

\noindent
{\bf V:}
$\int_1^\infty\, \frac{ \vartheta ( x) - x }{ x^2}\, dx$ est
une intégrale convergente.

\noindent
{\bf VI:} $\vartheta ( x) \sim x$.

\noindent
Le théorème des nombres premiers découle alors aisément de {\bf VI}, puisque, 
pour tout $\epsilon > 0$:
\[
\vartheta(x)
=
\sum_{p\leqslant x}\,{\rm log}\,p
\leqslant
\sum_{p\leqslant x}\,{\rm log}\,x
=
\pi(x)\,{\rm log}\,x
\]
\[
\aligned
\vartheta(x)
\geqslant
\sum_{x^{1-\epsilon}\leqslant p\leqslant x}
&
\geqslant
\sum_{x^{1-\epsilon}\leqslant p\leqslant x}\,
(1-\epsilon)\,{\rm log}\,x
\\
&
=
(1-\epsilon)\,{\rm log}\,x\,
\big[
\pi(x)+{\rm O}(x^{1-\epsilon})
\big].
\endaligned
\]

} en 1980, en modifiant astucieusement
les contours d'intégration de Hadamard et de la Vallée Poussin.
Par ailleurs,
en 1949, Selberg et Erd\"os ont élaboré une deuxième preuve
<<\,épurée\,>> qui évite le recours à l'analyse complexe, mais cette
preuve est relativement longue (une trentaine de pages) et elle ne
semble pas motiver ou offrir des développements ultérieurs.

\subsection*{ Démultiplication artificielle des énoncés}
Les démonstrations sont mobiles, transitoires, modifiables,
améliorables. Wittgenstein dit que toute nouvelle démonstration {\it
fabrique}\, une nouvelle connexion. Rien de plus exact. Mais dans
certains cas de figure tels que le théorème des nombres premiers, la
position wittgensteinienne absolutiste s'expose à une difficulté
spéculative qui lui sera immédiatement objectée par tout mathématicien
en acte: comment maintenir l'irréductibilité de nature entre énoncé et
démonstration, quand l'énoncé en question dont on cherche une
démonstration nouvelle a déjà été démontré par plusieurs voies
rigoureuses? L'énoncé reste-t-il irréductiblement ouvert et
conjectural? Doit-on exiger de la philosophie des
mathématiques qu'elle respecte le principe logique
de non-contradiction?

Parfois, deux démonstrations distinctes fournissent deux théorèmes
essentiellement équivalents mais qui sont légèrement différents, leur
différence pouvant être exprimée visiblement dans les énoncés: raison
est alors donnée à Wittgenstein. C'est notamment le cas dans les
mathématiques contemporaines, fortes d'un extrême raffinement, où des
équipes en compétition internationales développent des approches
concurrentes et bien distinctes pour étudier un même type de
problèmes: la différence des techniques utilisées remonte alors
jusqu'aux énoncés dans les publications. Est nouveau tout résultat
dont la démonstration est nouvelle.

Mais pour maintenir la cohérence globale de sa posture philosophique,
Wittgenstein semble prétendre que deux énoncés sont réellement
distincts dès lors que leur démonstrations diffèrent. Mais que dire
des énoncés tels que $\pi ( n) \sim \frac{ n}{ {\rm log}\, n}$ qui
sont exactement les mêmes parce qu'ils sont déjà atomiques et simples?
Faut-il chercher à faire transparaître à tout prix les différences
argumentatives des preuves dans les énoncés? Comment penser les degrés
de la différence? Faut-il supprimer les énoncés et ne mémoriser que
les démonstrations spécifiques?

On se trouve ainsi ramené à une vaste question: qu'est-ce qu'un énoncé
(une proposition, un théorème) mathématique? Tous les mathématiciens
se posent la question suivante: quelle forme donner à un énoncé que l'on
vient de démontrer? Aucune réponse définitive ou dogmatique ne peut
être proposée. \'Elasticité stylistique et souplesse du langage
complexifient encore le jeu de la publication. Pensée et écriture
mathématiques sont incapables de fixer définitivement leur rhétorique.

Par convention au moins, l'énoncé doit extraire une information
synthétique essentielle. La règle usuelle veut que l'énoncé soit
relativement court par rapport à la démonstration, ce qui est la
plupart du temps le cas.

Mais en vérité, nous retrouvons ici un des caractères distinctifs
fondamentaux de l'irréversible-synthétique: c'est bien parce que les
mathématiques sont faites d'obstacles, c'est bien parce que les
problèmes à résoudre exigent d'escalader ou de contourner lentement
des montagnes que l'irréversible-synthétique existe et se divise en
énoncés et démonstrations. Il y a là encore un problème crucial et
très difficile que la philosophie des mathématiques ne doit pas avoir
la tentation d'occulter: comment l'irréversible-synthétique est-il
possible?

\subsection*{ Arguments heuristiques en théorie analytique des nombres}
En suivant Hardy\footnote{\, {\em Ramanujan. Twelve lectures on
subjects suggested by his life and work}, Chelsea, New York, 1940.},
restituons deux arguments heuristiques simples qui conduisent à
l'équivalence $\pi ( n) \sim \frac{ n}{{\rm log}\, n}$ du théorème des
nombres premiers, ou, ce qui revient au même, à la conclusion:
\[
p_m\sim
m\,{\rm log}\,m,
\]
où $p_m$ est le $m$-ième nombre premier.

Voici le premier argument. Partons de l'identité d'Euler, 
valable uniformément pour ${\rm Re}\, s >1$:
\[
\small
\prod_p\,\frac{1}{1-p^{-s}}
=
\frac{1}{(1-2^{-s})(1-3^{-s})(1-5^{-s})\cdots}
=
\frac{1}{1^s}
+
\frac{1}{2^s}
+
\frac{1}{3^s}
+
\cdots
=
\sum_{m\geqslant 1}\,\frac{1}{m^s},
\]
où le produit porte sur l'ensemble des nombres premiers $\mathcal{ P }
= \{ 2, 3, 5, 7, 11, \dots \}$.  Il est naturel que le
produit $\prod \, \frac{ 1}{ 1 - p^{ -s}}$ et la série $\sum\, \frac{
1}{ m^s}$ divergent de la même manière\footnote{\, La manière dont ces
deux quantités divergent est nécessairement la même, puisque
l'identité d'Euler est valide quel que soit $s$ satisfaisant ${\rm
Re}\, s > 1$. Toutefois, c'est en estimant la contribution principale
de divergence pour chacun des deux membres qu'on peut être amené soit
à commettre une erreur, soit (si on ne s'est pas trompé sur le choix
des termes divergents principaux, ce qui est le cas ici) à éprouver
de réelles difficultés à transformer le raisonnement heuristique
en démonstration rigoureuse. }
lorsque $s$ tend vers $1$ en restant dans le demi-plan ${\rm Re} s >
1\}$.  Clairement, la série tronquée $\sum_{ m \leqslant
n}\, \frac{ 1}{ m}$ diverge comme ${\rm log}\, n$.  Par ailleurs, si
on développe le logarithme du produit:
\[
\small
\aligned
{\rm log}\,\prod_p\,
\frac{1}{1-p^{-s}}
&
=
\sum_p\,{\rm log}\,
\frac{1}{1-p^{-s}}
\\
&
=
\sum_p\,\frac{1}{p^s}
+
\bigg(
\sum_p\,\frac{1}{2\,p^{2s}}
+
\sum_p\,\frac{1}{3\,p^{3s}}
+
\cdots
\bigg)
\endaligned
\]
en tenant compte du fait que tous les termes $\sum_p\, \frac{ 1}{ k\,
p^{ ks}}$ pour $k\geqslant 2$ convergent, on s'attend à ce que,
lorsque $s \to 1$, la première somme $\sum\, \frac{ 1}{ p}$ diverge
comme ${\rm log}\, \Big( \sum\, \frac{ 1}{ m} \Big)$, ou plus
précisément:
\[
\small
\sum_{p\leqslant n}\,
\frac{1}{p}
\sim
{\rm log}\,
\bigg(
\sum_{m\leqslant n}\,
\frac{1}{m}
\bigg)
\sim
{\rm log}\,{\rm log}\,n.
\]
Comme par ailleurs:
\[
\sum_{m\leqslant n}\,\frac{1}{m\,{\rm log}\,m}
\sim
{\rm log}\,{\rm log}\,n,
\]
cette dernière formule pourrait indiquer que $p_m$ est
asymptotiquement égal à 
$m\, {\rm log}\, m$, ce qu'on voulait obtenir.

\subsection*{ Métaphysique des raisonnements heuristiques}
Ici, pétition de principe: sachant que la série de Bertrand $\sum_{
m\leqslant n}\, \frac{ 1}{ m\, {\rm log}\, m}$ diverge comme ${\rm
log}\, {\rm log}\, n$, on annonce que le comportement asymptotique
{\it inconnu}\, de $\frac{ 1}{ p_m}$ doit être le même que $\frac{ 1}{
m\, {\rm log}\, m}$.  Mais il se pourrait très bien qu'une infinité
d'autres séries différentes $\sum_m \, \frac{ 1}{ q_m}$ de termes
généraux $\frac{ 1}{ q_m}$ essentiellement distincts de $\frac{ 1}{
m\, {\rm log}\, m}$ diverge aussi comme ${\rm log}\, {\rm log}\, n$.
Dans l'absolu, ce raisonnement très périlleux devrait donc être
considéré comme irrecevable, à cause de la diversité {\it a priori}\,
du possible: les grandes catégories métaphysiques 
restent omniprésentes en mathématiques.

Mais à l'époque où écrit Hardy, la loi attendue $\pi ( n) \sim \frac{
n}{ {\rm log}\, n}$ (ou, de manière équivalente, $p_m \sim m \, {\rm
log}\, m$) avait été anticipée sur le plan expérimental depuis plus
d'un siècle par Legendre, Gauss et d'autres\,\,---\,\,sans compter que
les démonstrations rigoureuses de Hadamard et de la Vallée Poussin
circulaient depuis plus d'une vingtaine d'années.  On assiste donc ici
à un phénomène intéressant de {\it cristallisation de la cohérence}.
La spéculation mathématique est une machine à voyager dans le temps
irréversible du démonstratif.  Elle scrute
librement l'embryogénèse du déductif.

\section*{ Conclusion ouverte}

\subsection*{ \'Epilogue critique}
Le philosophe analytique wittgensteinien n'a peut-être pas encore pris
conscience du fait que la distinction capitale entre conjectures et
démonstrations n'est guère qu'une pièce initiale dans un puzzle
mathématique indéfini. Se crisper sur cette distinction expose aux
circularités spéculatives, aux répétitions désorganisées.  En
mathématiques, parce que tout se ramifie au-delà des racine, on
structure la pensée (c'est une règle d'or), on élimine l'extrinsèque,
on désigne l'inconnu, et on travaille au front. Le mathématicien en
acte joue en permanence avec les grands concepts de la métaphysique
classique: {\it a priori}\big/{\it a posteriori}; jugement
analytique\big/jugement synthétique; irréversible-synthétique;
dialectique; heuristique.  Et ses pensées jouent avec souplesse du
technique comme du médidatif.

\subsection*{ Penser le calcul}
Depuis une décennie, les logiciels de calcul formel tels que Maple,
Mathematica, Singular, Macaulay, Pari, {\it etc.}  sont régulièrement
enseignés dans les cursus de la Licence. Beaucoup de démonstrations
sont maintenant assistées par ordinateur. La recherche s'hybride.  En
géométrie, la pensée du continu ramifie ses discrétisations
conceptuelles. Tous ces éléments ne sont pas encore pensée par la
philosophie comme ils devraient l'être.

\subsection*{ Directions ouvertes de philosophie des mathématiques} \

\begin{itemize}

\smallskip\item[$\bullet$]
\'Edifier une \green{\sl pensée de l'ouverture mathématique technique}.

\smallskip\item[$\bullet$]
Spéculer sur la nature des questions mathématiques.

\smallskip\item[$\bullet$]
Ramifier la question kantienne: <<\,comment les jugements synthétiques
{\it a priori}\, sont-ils possibes\,>>.

\smallskip\item[$\bullet$]
Penser l'irréversible-synthétique.

\smallskip\item[$\bullet$]
Actifier, reproduire, propager, mécaniser, automatiser et désacraliser
le questionnement mathématique.

\smallskip\item[$\bullet$]
Formuler expressément les ouvertures rémanentes.

\smallskip\item[$\bullet$]
Constituer des catégories de pensée pour systématiser la nature de ce
qui demeure dans le domaine du non-exploré.

\smallskip\item[$\bullet$] 
Désigner l'indécision.

\smallskip\item[$\bullet$] 
Typiser et hiérarchiser les questions spécifiques.

\smallskip\item[$\bullet$] 
Démasquer les ignorances paradoxales qui se présentent com\-me
connaissances entrevues.

\smallskip\item[$\bullet$]
Réhabiliter le philosophique des mathématiques.

\end{itemize}

\vfill\end{document}